\documentclass[11pt]{amsart}
\usepackage[dvipdfmx]{graphicx}
\usepackage{amssymb,amstext,amsfonts,amscd}
\usepackage{amsmath}
\usepackage{latexsym}
\usepackage{amsthm}
\usepackage{multicol}
\usepackage{fancybox}
\usepackage{here}
\usepackage{multirow}
\usepackage{mathrsfs}

\def\opn#1#2{\def#1{\operatorname{#2}}}
\opn\Max{Max}
\opn\max{max}
\opn\Min{Min}
\opn\min{min}
\opn\rank{rank}
\opn\Cone{Cone}
\opn\Int{Int} %%% new
\opn\div{div}
\opn\Div{Div}
\opn\Ker{Ker}
\opn\id{id}
\opn\GL{GL}
\opn\SL{SL}
\opn\mod{mod}
\opn\det{det}
\opn\grad{grad}
\opn\codim{codim}
\opn\sign{sign}

\newcommand{\bz}{\mathbb{Z}}

\newcommand{\br}{\mathbb{R}}
\newcommand{\bc}{\mathbb{C}}

\newcommand{\zero}{\boldsymbol{0}}
\newcommand{\ba}{\boldsymbol{a}}   %%%% boldsymbol

\newcommand{\bbu}{\boldsymbol{u}}
\newcommand{\bv}{\boldsymbol{v}}
\newcommand{\bw}{\boldsymbol{w}}
\newcommand{\bx}{\boldsymbol{x}}
\newcommand{\by}{\boldsymbol{y}}
\newcommand{\bbz}{\boldsymbol{z}}

\newcommand{\La}{\mathcal{L}}

\newcommand{\Va}{\mathcal{V}}

\newcommand\si{\sigma}

\newcommand\De{\Delta}

\newcommand\rdeg{{\rm{rdeg}\/}}
\newcommand\pdeg{{\rm{pdeg}\/}}

\newcommand\inv{^{-1}}

\def\inv{^{-1}}

\newtheorem{theorem}{Theorem}
\newtheorem{definition}[theorem]{Definition}
\newtheorem{corollary}[theorem]{Corollary}
\newtheorem{proposition}[theorem]{Proposition}
\newtheorem{lemma}[theorem]{Lemma}
\newtheorem{remark}[theorem]{Remark}

\numberwithin{equation}{section}

\begin{document}
\pagestyle{plain}

\title
[RESOLUTIONS OF NEWTON NON-DEGENERATE MIXED POLYNOMIALS]   %%%
{RESOLUTIONS OF NEWTON NON-DEGENERATE MIXED POLYNOMIALS \\ %%%
OF STRONGLY POLAR NON-NEGATIVE \\                          %%%
MIXED WEIGHTED HOMOGENEOUS FACE TYPE}

\author{Sachiko Saito and Kosei Takashimizu}

\address{Department of Mathematics Education, Asahikawa Campus, Hokkaido University of Education, 
Asahikawa 070-8621, Japan}

\email {saito.sachiko@a.hokkyodai.ac.jp}


\begin{abstract}
Let $f(\bbz,\bar\bbz)$ be 
a convenient Newton non-degenerate mixed polynomial with strongly polar non-negative mixed weighted homogeneous face functions. 
We consider a convenient regular simplicial cone subdivision $\Sigma^*$ which is admissible for $f$ 
and take the toric modification $\hat{\pi} : X \to \bc^n$ associated with $\Sigma^*$. 
We show that the toric modification resolves topologically 
the singularity of the mixed hypersurface germ defined by $f(\bbz,\bar\bbz)$ 
under the {\sc Assumption}{\rm (*)} (Theorem \ref{theorem11-improved}). 
This result is an extension of the first part of Theorem 11 (\cite{Oka2015}) by M. Oka, 
which studies strongly polar positive cases, 
to strongly polar non-negative cases. 
We also consider some typical examples (\S \ref{section-exam4-3-mixed}). 
\end{abstract}

\maketitle

%%%%%%%%%%%%%%%%%%%%%%%%%%%%%%%%%%%%%%%%%%%%%%%%%%%%

\section{Introduction}\label{section-intro}

A {\em mixed analytic function} (or {\em mixed function}) $f(\bbz, \bar{\bbz})$ on 
a neighborhood $U$ of $\zero$ in $\bc^n$ is defined to be 
$f(\bbz, \bar{\bbz}) := F(\bbz, \bar{\bbz})$, where 
$F(\bbz, \bw)$ on $U\times U$ is a complex valued holomorphic function with complex $2n$ variables. 
For a mixed function $f(\bbz, \bar{\bbz})$ on $U$, we consider the {\em mixed hypersurface} 
$V := f^{-1}(0)$ as a germ at $\zero$. 
We study resolutions of the singularity $(V,\zero)$ by {\em toric modifications} 
when $\zero$ is an isolated {\em mixed singular point {\rm (see \S \ref{section-mixed-function})} of $V$}. 
In \cite{Oka2015}, M. Oka obtained Theorem 11 concerning 
topological resolutions of the isolated singularity $\zero \in V := f^{-1}(0)$ 
by toric modifications under the assumption that $(f,\zero)$ is 
a {\em convenient} and {\em Newton non-degenerate} mixed polynomial germ 
{\em of strongly polar \underline{positive} mixed weighted homogeneous face type}. 
Theorem 11 of \cite{Oka2015} (see also Theorem 9.19 of \cite{Oka2018}) 
is a generalization of Varchenko's famous Lemma 10.3 (\cite{Varchenko1976}) on holomorphic function germs 
to mixed polynomial germs. 

The main purpose of this paper is to prove an extension (Theorem \ref{theorem11-improved}) 
of the first part ((1) and (2)) of Theorem 11 of \cite{Oka2015} 
to {\em strongly polar \underline{non-negative} mixed weighted homogeneous face type} cases. 
Moreover, our Theorem \ref{theorem11-improved} includes the information about 
the differentiability ($C^s$-class, where $s\geq 0$) of the strict transform $\tilde V$ of $V$. 
At present we do not get any extensions of the latter part (3) of Theorem 11 (\cite{Oka2015}), 
which is a result on the zeta functions of the Milnor fibrations of mixed polynomials. 
In strongly polar non-negative mixed weighted homogeneous face type cases, 
we do not always have the Milnor fibrations. 
The simplest example is the mixed function $f : \bc \to \bc,\ f(z, \bar{z}):= z\bar{z}$. 
This function does not have the Milnor fibration but 
the zero set $f^{-1}(0)$ still coincides with the holomorphic case: $g(z):= z^2$. 
In \S \ref{section-exam4-3-mixed} of this paper we consider and calculate some typical examples including 
the mixed polynomial $(z_1^4 + z_2^3)\overline{(z_1^2 + z_2^3)}$. 
Before the studies of mixed functions, a real function of type $f\overline{g}$ had been studied by 
A. Pichon and J. Seade \cite{Pichon-Seade2008}. 
They mentioned that it has the Milnor fibration under the assumption that 
the multiplicities of $f$ and $g$ do not coincide on each of the exceptional divisors. 
The coincidence of the multiplicities means that 
the polar degree of the corresponding face function of $f\overline{g}$ is $0$. 
See also \cite{Pichon2005} and \cite{Pichon-Seade2011}. 

In \S \ref{section-mixed-function} and \S \ref{section-Newton polyhedron} 
we introduce some basic terms, namely, {\em mixed functions}, {\em mixed polynomials}, {\em mixed critical points}, 
{\em the (radial) Newton polyhedron} of a mixed function germ, the {\em face functions} of a mixed function germ, 
and so on. 
In \S \ref{section-whp} we give the definitions of 
mixed function germs {\em of strongly polar positive} (respectively, {\em of strongly polar non-negative}) {\em mixed weighted homogeneous face type}. 
In \S \ref{section-Newton-non-degeneracy} 
we describe the notion of {\em Newton non-degeneracy} of mixed function germs. 
We also mention about {\em strong Newton non-degeneracy}. 
Although we do not use strong Newton non-degeneracy 
in the statement of our main Theorem \ref{theorem11-improved}, 
we should note that 
for mixed weighted homogeneous polynomials (Definition \ref{mixed-whp}) of non-zero polar degree, 
Newton non-degeneracy implies strong Newton non-degeneracy. 
See Remark \ref{remark4}. 
In \S \ref{section-cone-subdivision} and \S \ref{section-convenient-function} 
we moreover introduce important terms, {\em regular simplicial cone subdivisions}, 
{\em toric modifications}, and {\em convenient} mixed function germs. 
We prepare the key Propositions \ref{non-empty}---\ref{intersection}, 
Lemma \ref{E-j} and Proposition \ref{change-of-strictly-positive-weight-vector}, 
which we use in the proof of the main Theorem \ref{theorem11-improved}. 

In \S \ref{section-main} we state the main Theorem \ref{theorem11-improved}. 
The assertion {\bf (ii)} of Theorem \ref{theorem11-improved} includes 
the differentiability of the strict transform germ $\tilde V$. 
For the proof of {\bf (ii)}, 
we prepare Lemma \ref{fraction_lemma}, 
which describes the differentiability of functions in the form of $u^{r+s} / \bar{u}^r$. 

In the statement of Theorem \ref{theorem11-improved} {\bf (i)}, 
we put the assumption that ``the vertices of $\tau$ are all strictly positive". 
We at present cannot prove {\bf (i)} without this assumption. 
Following {\bf (i)}, we put {\sc Assumption}{\rm (*)} in the statement of Theorem \ref{theorem11-improved} {\bf (ii)}. 
The reason is as follows: 
We prove {\bf (ii)} by using the argument in the proof of {\bf (i)}. 
For an arbitrary point $\hat{\bbu}^0 \in \tilde V \cap (\bigcup_{P\gg 0}\hat{E}(P))$, 
there exist some vertices $P_1, \dots , P_\kappa \ (1 \leq \kappa \leq n)$ and 
an $n$-dimensional cone $\sigma = (P_1, \dots , P_n)$ for which 
the toric chart $\bc_\sigma^n$ contains 
$\hat{\bbu}^0 = (\underbrace{0, \dots, 0}_{\kappa} \, ,\ u^0_{\kappa +1}, \dots, u^0_n)$ and 
$u^0_{\kappa +1}\neq 0, \dots ,u^0_n\neq 0$. 
If some of $P_1, \dots , P_\kappa$ are not strictly positive, 
then we cannot apply {\bf (i)} to the cone $(P_1, \dots , P_\kappa)$. %%%%
Therefore we assume {\sc Assumption}{\rm (*)} in {\bf (ii)}. 
It will be interesting to study whether {\sc Assumption}{\rm (*)} can be removed or not. 

In the last section \ref{section-exam4-3-mixed} we investigate 
the mixed polynomials $z_1^4 \bar{z}_1^2 + z_1^a \bar{z}_1^{4-a} z_2^b \bar{z}_2^{3-b} + \bar{z}_1^2 z_2^3 + z_2^3 \bar{z}_2^3$ 
including $(z_1^4 + z_2^3)\overline{(z_1^2 + z_2^3)}$ in detail. 

\medskip

\section{Mixed analytic functions and mixed polynomials}\label{section-mixed-function}

Let $U$ be a neighborhood of $\zero$ in $\bc^n$. 
We assume $\bar{U}=U$, 
where $\bar{\bbz}$ stands for 
the complex conjugate $(\bar{z_1}, \dots , \bar{z_n})$ of $\bbz=(z_1, \dots , z_n) \in \bc^n$. 
For a complex valued holomorphic function $F(\bbz, \bw)$ on $U\times U$ with complex $2n$ variables, 
we set 
$$
f(\bbz, \bar{\bbz}) := F(\bbz, \bar{\bbz}), 
$$
which is defined over $U$. 
We call it a {\em mixed analytic function} (or {\em mixed function}) $f$ on $U$. 
We assume that 
$
F(\zero, \zero) = 0
$
throughout this paper, 
and hence, 
$$f(\zero)=0.$$
Let 
$$
F(\bbz, \bw) = \sum_{\nu, \mu}c_{\nu, \mu} \bbz^\nu \bw^\mu
$$
be the Taylor expansion of $F$ at $(\zero , \zero)$, 
where 
$
\nu = (\nu_1, \dots , \nu_n),\ \mu = (\mu_1, \dots , \mu_n),\ 
\nu_i \geq 0,\ \mu_j \geq 0,\ \ 
\bbz^\nu := z_1^{\nu_1}\cdots z_n^{\nu_n},\ \ \bw^\mu := w_1^{\mu_1}\cdots w_n^{\mu_n} .
$
Then we have $c_{\zero, \zero} = 0$ and 
$$
f(\bbz, \bar{\bbz}) = \sum_{\nu, \mu}c_{\nu, \mu} \bbz^\nu \bar{\bbz}^\mu .
$$
We call $f(\bbz, \bar{\bbz})$ a {\em mixed polynomial} 
when the number of monomials $c_{\nu, \mu} \bbz^\nu \bar{\bbz}^\mu,\ c_{\nu, \mu}\neq 0$ is finite. 

For a mixed function $f$ on $U$, let us consider 
$$
V := f^{-1}(0) \ \ (\subset U), 
$$
which we call a {\em mixed hypersurface}.  
We set $z_j = x_j + iy_j$ and consider the real valued functions 
$$g(\bx, \by) := \Re f(\bbz, \bar{\bbz}),\ \ h(\bx, \by) := \Im f(\bbz, \bar{\bbz})$$
with real $2n$ variables $(\bx, \by)$. 
Then we have 
$
f(\bbz, \bar{\bbz}) = g(\bx, \by) + ih(\bx, \by), 
$
and 
$$
V=\{ (\bx, \by) \in U \ |\ g(\bx, \by)=h(\bx, \by)=0 \} .
$$
Thus $V$ is a real analytic variety in $U \ (\subset \br^{2n})$. 

\medskip

\begin{definition}\label{mixed-critical}
{\rm 
Let $f(\bbz, \bar{\bbz})$ be a mixed analytic function on $U \ (\subset \bc^n)$. 
We say 
$
\ba = (a_1, \dots , a_n) \in U
$
is a {\em mixed critical point} (or a {\em mixed singular point}) of $f$ 
if 
the rank of the differential map 
$(df)_{\ba} : T_{\ba} \bc^n \to T_{f(\ba)} \bc \cong T_{f(\ba)} \br^2$ 
is less than $2$. 
We also say 
$\ba$ is a {\em mixed singular point of the mixed hypersurface $f^{-1}(f(\ba, \bar{\ba}))$} 
when 
$\ba$ is a mixed critical point of $f$. 
We say $\ba \in U$ is a {\em mixed regular point} of $f$ 
if it is not a mixed critical point of $f$.
}
\end{definition}

\medskip

We set 
$
\partial f := \left(\frac{\partial f}{\partial z_1}, \dots , \frac{\partial f}{\partial z_n} \right),\ \ \ 
\bar{\partial} f := \left(\frac{\partial f}{\partial \bar{z_1}}, \dots , \frac{\partial f}{\partial \bar{z_n}} \right) ,
$
where 
\begin{equation}\label{partial1}
\frac{\partial}{\partial z_j} = 
\frac{1}{2}\left( \frac{\partial}{\partial x_j} - i \frac{\partial}{\partial y_j} \right),\ \ \ 
\frac{\partial}{\partial \bar{z_j}} = 
\frac{1}{2}\left( \frac{\partial}{\partial x_j} + i \frac{\partial}{\partial y_j} \right) .
\end{equation}

Equivalently, we have
\begin{equation}\label{partial2}
\frac{\partial}{\partial x_j} = 
\frac{\partial}{\partial z_j} + \frac{\partial}{\partial \bar{z_j}}\,,\ \ \ 
\frac{\partial}{\partial y_j} = 
i \left( \frac{\partial}{\partial z_j} - \frac{\partial}{\partial \bar{z_j}} \right) .
\end{equation}

\begin{proposition}[Oka \cite{Oka2008}, Proposition 1, \cite{Oka2018}, Proposition 8.1]\label{critical}
The following two conditions are equivalent:
\begin{enumerate}
\item $\ba = (a_1, \dots , a_n)$ is a mixed critical point of $f$. 
\item There exists a complex number $\alpha$ with $|\alpha|=1$ which satisfies 
$\overline{\partial f(\ba, \bar{\ba})} = \alpha \bar{\partial}f(\ba, \bar{\ba})$. 
\end{enumerate}
\end{proposition}

\begin{corollary}\label{regular-criterion}
If 
$$
\left|\frac{\partial f}{\partial z_{j}}(\ba, \bar{\ba})\right| \neq 
\left|\frac{\partial f}{\partial \overline{z_{j}}}(\ba, \bar{\ba})\right|
$$
for some $j \ (1 \leq j \leq n)$, 
then $\ba$ is a mixed regular point of $f$. 
\end{corollary}

\medskip

\section{(Radial) Newton polyhedrons and their faces}\label{section-Newton polyhedron}

We first define 
$
\br_+^n := \{ (x_1, \dots , x_n)\ |\ x_i \geq 0\}     %%%
$. 
Let $(f,\zero)$ be the germ of a mixed function 
$$
f(\bbz, \bar{\bbz}) = \sum_{\nu, \mu}c_{\nu, \mu} \bbz^\nu \bar{\bbz}^\mu
$$
at $\zero \in \bc^n$. 
Let 
$$
\Gamma_+(f)
$$
be the convex hull of the set 
$$
\displaystyle \bigcup_{c_{\nu, \mu} \neq 0}  (\nu + \mu) + \br_+^n ,
$$
which we call 
the {\em (radial) Newton polyhedron} of the germ $(f,\zero)$ of a mixed function $f$ at $\zero$. 

\begin{definition}[face of $\Gamma_+(f)$]\label{face-mixed}
{\rm 
For a ``weight vector" 
$$
P={}^t(p_1, \dots , p_n)\ (\neq \zero) \in N^+_\br := \br_+^n,    %%%
$$
let 
$$d(P)$$
be the minimum value of the linear function 
$$
P : \Gamma_+(f) \to \br,\ P(\xi) := \sum_{j=1}^n p_j \xi_j \ \ (\in \br),
$$
where 
$\xi = (\xi_1 , \dots , \xi_n) \in \Gamma_+(f)$. 
We set 
$$
\Delta(P) := \{ \xi \in \Gamma_+(f) \ |\ P(\xi) = d(P) \},    %%%
$$
which we call a {\em face} of $\Gamma_+(f)$. 
Note that $\Delta(P) \neq \emptyset$ by its definition. 
If $\dim \Delta(P) = n-1$, then $P$ is unique up to multiplications of positive real numbers. 
}
\end{definition}

\begin{definition}[strictly positive weight vector]
{\rm 
We say a weight vector $P={}^t(p_1, \dots , p_n)$ is {\em strictly positive} if %%
$p_i >0$ for every $i \ (=1, \dots , n)$. 
We write $P\gg 0$ if $P$ is strictly positive.
}
\end{definition}

\medskip

Note that a face $\Delta$ of $\Gamma_+(f)$ is compact 
if and only if 
$\Delta = \Delta(P)$ 
for some strictly positive weight vector $P$. For the proof, for example, see \cite{Saito-Takashimizu2021-8}.

\medskip

\begin{definition}[the (radial) Newton boundary $\Gamma(f)$]
{\rm 
Let $\Gamma(f)$ be the union of all compact faces of $\Gamma_+(f)$, 
which we call 
the {\em (radial) Newton boundary} of the germ $(f,\zero)$ of a mixed function $f$. 
}
\end{definition}

\begin{definition}[face function]\label{face-function}
{\rm 
For a \underline{compact} face $\Delta(P) \ (\subset \Gamma(f))$, we define 
$$
f_{\Delta(P)}(\bbz) \ (\text{or}\ f_P(\bbz)) := \sum_{\mu + \nu \in \Delta(P)}c_{\nu, \mu} \bbz^\nu \bar{\bbz}^\mu ,
$$
which we call a {\em face function ({\rm or} face polynomial)} of a mixed function germ $(f,\zero)$ (cf. \cite{Oka2018}, p.78). 
}
\end{definition}

\medskip

\section{Radial and polar weighted homogeneous polynomials}\label{section-whp}

\begin{definition}[\cite{Oka2018}, p.182, \cite{Oka2015}]
{\rm 
A mixed polynomial $f(\bbz, \bar{\bbz}) = \sum_{\nu, \mu}c_{\nu, \mu} \bbz^\nu \bar{\bbz}^\mu$ is called 
{\em radially weighted homogeneous} if 
there exists a weight vector 
$$Q={}^t(q_1, \dots , q_n)\ (\neq \zero) \in N^+_\br := \br_+^n$$
and 
a positive integer $d_r \ (> 0)$ such that 
$$
c_{\nu, \mu} \neq 0 \Longrightarrow \sum_{i=1}^n q_i(\nu_i + \mu_i) = d_r .
$$
We call $d_r$ {\em the radial degree} of $f$, and define 
$$
\rdeg_Q f := d_r .
$$

A mixed polynomial $f(\bbz, \bar{\bbz}) = \sum_{\nu, \mu}c_{\nu, \mu} \bbz^\nu \bar{\bbz}^\mu$ is called 
{\em polar weighted homogeneous} if 
there exists a weight vector 
$$P={}^t(p_1, \dots , p_n)\ (\neq \zero) \in N^+_\br := \br_+^n$$
and 
an integer $d_p$ ($>0,\ 0\ \text{or}\ <0$) such that 
$$
c_{\nu, \mu} \neq 0 \Longrightarrow \sum_{i=1}^n p_i(\nu_i - \mu_i) = d_p .
$$
We call $d_p$ {\em the polar degree} of $f$, and define 
$$
\pdeg_P f := d_p .
$$
}
\end{definition}

\begin{remark}
{\rm 
Since we assume $f(\zero)=0$, 
every face function (Definition \ref{face-function}) $f_P(\bbz) = f_{\Delta(P)}(\bbz)$, 
where we take a strictly positive $P$, 
of a mixed function germ $(f,\zero)$ 
is a radially weighted homogeneous polynomial of radial degree $d(P)\ (>0)$ 
with respect to the weight vector $P$. 
}
\end{remark}

\begin{definition}[\cite{Oka2018}, pp.182--184]\label{mixed-whp}
{\rm 
We say a mixed polynomial $f(\bbz, \bar{\bbz})$ is 
a {\em mixed weighted homogeneous polynomial}
\footnote{We use the definition in \cite{Oka2018}.}
if it is both radially and polar
\footnote{The corresponding weight vectors $Q$ and $P$ are possibly different.}
 weighted homogeneous.
}
\end{definition}

\begin{definition}[cf. \cite{Oka2018},p.183,\ \cite{Oka2015},p.174]\label{strongly-mixed-whp}
{\rm 
We say a mixed weighted homogeneous polynomial $f(\bbz, \bar{\bbz})$ is 
a {\em strongly mixed weighted homogeneous polynomial} if 
$f$ is radially and polar weighted homogeneous 
with respect to the \underline{same} weight vector $P$. 
Furthermore, $f$ is called 
a {\em strongly polar positive} (respectively, {\em strongly polar non-negative}) {\em mixed weighted homogeneous polynomial} 
if $\pdeg_P f >0\ (\text{respectively,}\ \pdeg_P f \geq 0)$. 
}
\end{definition}

\begin{definition}[cf. \cite{Oka2018},Definition 9.18,\ \cite{Oka2015},p.174]\label{strongly-mixed-wh-face-type}
{\rm 
Let $(f,\zero)$ be a mixed function germ at $\zero \in \bc^n$. 
\begin{enumerate}
\item The germ $(f,\zero)$ of a mixed function $f(\bbz, \bar{\bbz})$ at $\zero \in \bc^n$ is called 
{\em of strongly mixed weighted homogeneous face type} if 
the face function $f_\Delta(\bbz, \bar{\bbz})$ is 
a strongly mixed weighted homogeneous polynomial 
(Definition \ref{strongly-mixed-whp}) 
for every \underline{compact} face $\Delta$. 

\item The germ $(f,\zero)$ of a mixed function $f(\bbz, \bar{\bbz})$ at $\zero \in \bc^n$ is called 
{\em of strongly polar positive} (respectively, {\em of strongly polar non-negative}) {\em mixed weighted homogeneous face type} if 
for every 
%%%%%%%%%%%%%%%%%%%%%%%%%%%%%%%%%%%%%%%%%%%%%%%%
\footnote{In \cite{Oka2015}, a convenient (see Definition \ref{convenient-function} of this paper) mixed function $f(\bbz, \bar{\bbz})$ is called %%
{\em of strongly polar positive (mixed) weighted homogeneous face type} if 
the face function $f_\Delta(\bbz, \bar{\bbz})$ is a strongly polar positive (mixed) weighted homogeneous polynomial 
for every $(n-1)$-dimensional face. %%
After this definition, Proposition 10 of \cite{Oka2015} proves that 
for a convenient mixed function $f$ of strongly polar positive (mixed) weighted homogeneous face type (in the sense of \cite{Oka2015}) and 
any weight vector $P$, 
the face function $f_\Delta(P)$ is also a strongly polar positive mixed weighted homogeneous polynomial with respect to $P$.}
%%%%%%%%%%%%%%%%%%%%%%%%%%%%%%%%%%%%%%%%%%%%%%%%
compact face $\Delta$, 
the face function $f_\Delta(\bbz, \bar{\bbz})$ is 
a strongly polar positive (respectively, strongly polar non-negative) mixed weighted homogeneous polynomial 
(Definition \ref{strongly-mixed-whp}) 
with respect to 
some strictly positive weight vector $P$ with $\Delta = \Delta(P)$. %%
\end{enumerate}
}
\end{definition}

\begin{remark}
{\rm 
For mixed function germs of strongly polar non-negative mixed weighted homogeneous face type, 
Proposition \ref{change-of-strictly-positive-weight-vector} is important. 
}
\end{remark}

\medskip

\section{Newton non-degeneracy}\label{section-Newton-non-degeneracy}

\begin{definition}[\cite{Oka2010}, p.6, Definition 3,\ \cite{Oka2018}, p.80 and pp.181--182]
{\rm 
Let $(f,\zero)$ be the germ of a mixed function $f$ at $\zero \in \bc^n$. 
\begin{enumerate}
\item We say $(f,\zero)$ is {\em Newton non-degenerate 
over a compact face $\Delta \ (\subset \Gamma(f))$} 
if 
$0$ is not a mixed critical value of the face function 
$f_\Delta : {\bc^*}^n \to \bc$. 
(In particular, 
if $f_\Delta^{-1}(0) \cap {\bc^*}^n = \emptyset$, then $0$ is not a mixed critical value of the face function 
$f_\Delta : {\bc^*}^n \to \bc$.) 

\item Let $\Delta \ (\subset \Gamma(f))$ be a compact face with $\dim \Delta \geq 1$. 
We say $(f,\zero)$ is {\em strongly Newton non-degenerate 
over $\Delta$} 
if 
the face function $f_\Delta : {\bc^*}^n \to \bc$ has no mixed critical points 
and $f_\Delta : {\bc^*}^n \to \bc$ is surjective onto $\bc$. 

\item Let $\Delta \ (\subset \Gamma(f))$ be a compact face with $\dim \Delta = 0$, that is, 
$\Delta$ is a vertex of $\Gamma_+(f)$. 
We say $(f,\zero)$ is {\em strongly Newton non-degenerate 
over $\Delta$} 
if 
the face function $f_\Delta : {\bc^*}^n \to \bc$ has no mixed critical points
\footnote{In this case we do not need the surjectivity from ${\bc^*}^n$ onto $\bc$.}
. 
\end{enumerate}
}
\end{definition}

\begin{definition}[\cite{Oka2018}, p.80 and p.182]
{\rm 
We say the germ $(f,\zero)$ of a mixed function $f$ at $\zero \in \bc^n$ is 
{\em Newton non-degenerate} (respectively, {\em strongly Newton non-degenerate}) %%%
if 
$(f,\zero)$ is Newton non-degenerate (respectively, strongly Newton non-degenerate) 
over every compact face $\Delta \ (\subset \Gamma(f))$. 
}
\end{definition}

%%%%%%%%%%%%%%%%%%%%%%%%%%%%%%%%%%%%%%%%%%%%%%%%%%%%%%%%%%%%%%%%%%%%%%%%%%%%%%%%%%%%%%%%
\begin{remark}[cf. \cite{Oka2010}, Remark 4,\ \cite{Takashimizu2021thesis},\ \cite{Saito-Takashimizu2021-8}]\label{remark4}
{\rm 
Let $f(\bbz)$ be a holomorphic weighted homogeneous polynomial %%
of positive degree with respect to a strictly positive weight vector $P$. 
Then we have $f = f_{\Delta(P)}$.  %%%

(i)\ Suppose that $f$ is Newton non-degenerate over $\Delta(P)$, %%%
namely, 
$0$ is not a critical value of $f_{\Delta(P)}=f : {\bc^*}^n \to \bc$. 
Then we can show that 
$f_{\Delta(P)}=f : {\bc^*}^n \to \bc$ has no critical point. 
Hence, using (ii) below, we see that $f$ is strongly Newton non-degenerate over $\Delta(P)$. 

(ii)\ Suppose that $\dim \Delta(P) \geq 1$, namely, $f = f_{\Delta(P)}$ has at least two monomials. 
Then we can show that $f_{\Delta(P)}=f : {\bc^*}^n \to \bc$ is surjective. 

\medskip

In the mixed cases, the situations are more complicated. 
Let $(f,\zero)$ be the germ of a mixed function $f$ at $\zero \in \bc^n$. %%%

(iii)\ Suppose that $(f,\zero)$ is Newton non-degenerate over a compact face $\Delta = \Delta(P) \ (\subset \Gamma(f))$ and 
the face function $f_P$ is 
a mixed weighted homogeneous polynomial (Definition \ref{mixed-whp}). 
Morerover, we assume $P\gg 0$ and 
$f_P$ is a polar weighted homogeneous polynomial of \underline{non-zero} polar degree with respect to some weight vector $Q$. 
Then, we can show that $f_P : {\bc^*}^n \to \bc$ has no mixed critical points. 

(iv)\ In addition to (iii), %%%
we assume that 
$f_P^{-1}(0) \cap {\bc^*}^n \neq \emptyset$.  %%
Then we can show that $f_P : {\bc^*}^n \to \bc$ is surjective. 
Hence, with tha fact (iii), in this case, 
Newton non-degeneracy over a compact face $\Delta(P)$ implies strong Newton non-degeneracy over $\Delta(P)$. 

The proofs of (i) -- (iv) are given in \cite{Takashimizu2021thesis} and \cite{Saito-Takashimizu2021-8} in detail. 
}
\end{remark}

\medskip

\section{Admissible regular simplicial cone subdivisions}\label{section-cone-subdivision}

A {\em regular simplicial cone subdivision} of $N^+_\br$ is defined as in \cite{Oka2018}, p.86, \S 5.6. 
Let us introduce the notions of 
{\em convenient} subdivisions. 
We define 
$$
E_j := {}^t(0, \dots , 0, 1, 0, \dots , 0)
$$
for $j \ (= 1, \dots , n)$. 
Namely, $\{ E_1, \dots , E_n \}$ is the canonical basis of $\bz^n$. 

\begin{definition}[\cite{Oka2018}, p.90] \label{convenient-subdivision} 
{\rm 
We say a regular simplicial cone subdivision $\Sigma^*$ is {\em convenient} 
if 
$\Sigma^*$ contains all cones 
$$E_I := \Cone (E_{i_1}, \dots , E_{i_k})$$
with $I := \{ i_1, \dots , i_k \} \ (\subsetneqq \{1, \dots, n\})$. 
Thus $\Sigma^*$ contains $\Cone(E_1), \dots , \Cone(E_n)$. %%%
}
\end{definition}

\medskip

Let $\mathcal{V}$ be the set of all vertices of a regular simplicial cone subdivision $\Sigma^*$, and 
$\mathcal{V}^{+}$ be the set of all strictly positive vertices of $\Sigma^*$. 

\begin{remark} \label{convenient-subdivision-remark}
{\rm 
For a convenient regular simplicial cone subdivision $\Sigma^*$, the following facts are important:
\begin{enumerate}
\item The vertices except $E_1, \dots , E_n$ of $\Sigma^*$ are all strictly positive
\footnote{The proof is as follows. Let $P \ (\neq \zero)$ be a vertex in $\mathcal{V} \setminus \{ E_1, \dots , E_n \}$. 
Suppose that $P$ is not strictly positive. Then we may assume that $P = {}^t(p_1, \dots , p_n)$ with 
$p_1 >0, \dots , p_s >0$ and $p_{s+1} = \cdots = p_n = 0$. Here we have $2 \leq s < n$. 
Since $\Sigma^*$ is convenient and $s < n$, we see $\Cone (E_1, \dots , E_s) \in \Sigma^*$. 
We also see $\Cone (P) \neq \Cone (E_1, \dots , E_s)$ as $s\geq 2$. 
Then, by the definition of simplicial cone subdivisions (see \cite{Oka2018}, p.86), %%
we have $P \in \Cone (P) \cap \Cone (E_1, \dots , E_s) \subset \partial (\Cone (E_1, \dots , E_s))$. 
This is a contradiction. Hence, $P$ is strictly positive.}
. 
Namely, $\mathcal{V} \setminus \{ E_1, \dots , E_n \} \subset \mathcal{V}^{+}$. %% 
Thus we have %%
$$\mathcal{V} = \mathcal{V}^{+} \cup \{ E_1, \dots , E_n \}.$$

\item The toric modification (see \cite{Oka1997},p.73, \cite{Oka2010} or \cite{Oka2018},p.89 for the definition) 
$$
\hat{\pi} : X \to \bc^n
$$
associated with $\Sigma^*$ 
is a proper birational map and 
$
\hat{\pi} : X\setminus \hat{\pi}^{-1}(\zero) \to \bc^n \setminus \{ \zero \}
$
is biholomorphic and we have 
$$
\hat{\pi}^{-1}(\zero) = \bigcup_{P\gg 0}\hat{E}(P),
$$
where $\hat{E}(P)$ is the exceptional divisor (\cite{Oka1997}, p.74) defined by a vertex $P$ of $\Sigma^*$. 
See \cite{Oka1997},Theorem(1.4),\ Corollary(1.4.1), and also \cite{Oka2018},Theorem 5.16. 
\end{enumerate}
}
\end{remark}

%%%%%%%%%%%%%%%%%%%%%%%%%%%%%%%%%%%%%%%%%%%%%%%%%%%%%%%%%%%%%%%%%%%%%%

\bigskip

Now for the germ $(f,\zero)$ of a mixed function $f$, we can also define 
the {\em (radial) dual Newton diagram} 
$\Gamma^*(f)$ 
of $(f,\zero)$ in the same way as holomorphic function germs at $\zero$
(see \cite{Oka1997},p.120,\ \cite{Oka2018}, p.97). 

\begin{definition}[admissible cone,\ \cite{Oka1997},p.121]
{\rm 
Let $\tau = \Cone(P_1, \dots , P_s)$, where $P_1, \dots , P_s \in N^+_\br$, %%%%
be a cone in $N^+_\br$, that is, 
$$
\tau := \{ \sum_{j=1}^s r_j P_j \ |\ r_j \geq 0,\ j = 1, \dots ,s \} .
$$
We set 
$$
\Int (\tau) := \{ \sum_{j=1}^s r_j P_j \ |\ r_1 >0, \dots ,r_s >0 \} ,
$$
which we call the {\em interior} of $\tau$, 
and 
$$
P_\tau := P_1 + \cdots + P_s .  %%%
$$
We say (cf. \cite{Oka1997},p.121,Definition (3.1.3)) 
a cone $\tau$ is {\em admissible} for $\Gamma^*(f)$ if 
the interior $\Int (\tau)$ %%%
is a subset of an equivalence class in $\Gamma^*(f)$, that is, 
$\Int (\tau)$ is contained in the equivalence class $[P_\tau]$. %%%
}
\end{definition}

\begin{definition}[admissible regular simplicial cone subdivision,\ \cite{Oka1997},p.121,\ \cite{Oka2018},p.98]\label{admissible-regul-simpl-cone-subdivision}
{\rm 
Let $\Sigma^*$ be a regular simplicial cone subdivision of $N^+_\br$. 
We say that $\Sigma^*$ is {\em admissible} for the germ $(f,\zero)$ of a mixed function $f$ (or for $\Gamma^*(f)$) 
if 
every cone $\tau$ in $\Sigma^*$ is admissible for $\Gamma^*(f)$. 
}
\end{definition}

\begin{proposition}[\cite{Oka1997}, p.121, Proposition (3.1.4)]\label{non-empty}
Let $\tau = \Cone(P_1, \dots , P_s)$, where $P_1, \dots , P_s \in N^+_\br$, %%%%
be a cone in $N^+_\br$. 
If $\tau$ is admissible for $\Gamma^*(f)$, 
then 
$
\De(P_i) \supset \De(P_\tau)
$
for any $i = 1, \dots , s$, where $P_\tau := P_1 + \cdots + P_s$. 
In particular, we have
$$
\bigcap_{i=1}^s \De(P_i) \supset \De(P_\tau) = \De(P) \ (\neq \emptyset)
$$
for every $P$ in $\Int (\tau)$, 
and therefore we have 
$
\bigcap_{i=1}^s \De(P_i) \neq \emptyset    %%%
$. 
\end{proposition}

\begin{proposition}\label{intersection-of-faces}
Let $\tau = \Cone(P_1, \dots , P_s)$, where $P_1, \dots , P_s \in N^+_\br$, %%%%
be a cone in $N^+_\br$. 
If $\tau$ is admissible for $\Gamma^*(f)$, 
then we have 
$$
\bigcap_{i=1}^s \De(P_i) = \De(P)
$$
for every $P$ in $\Int (\tau)$. 
\end{proposition}

\begin{proof}
By Proposition \ref{non-empty}, we have 
$\bigcap_{i=1}^s \De(P_i) \supset \De(P_\tau) \neq \emptyset$, where $P_\tau := P_1 + \cdots + P_s$. %%
Hence, take a point $\bx \in \De(P_\tau)$. 
Then, since $\bx \in \De(P_i)$ for any $i = 1, \dots , s$, 
we have 
$$
P_\tau(\bx)  =  \sum_{i=1}^s P_i(\bx)
             =  \sum_{i=1}^s d(P_i) .
$$
On the other hand, for any $\by \in \Gamma_+(f)$, we have 
$$
P_\tau(\by) = \sum_{i=1}^s P_i(\by) \geq \sum_{i=1}^s d(P_i) .
$$
Thus we conclude that 
$$
d(P_\tau) = \sum_{i=1}^s d(P_i) .  %%
$$

Let $\bx'$ be an arbitrary point of $\bigcap_{i=1}^s \De(P_i)$. 
Then we have 
$$
P_\tau(\bx') = \sum_{i=1}^s P_i(\bx') = \sum_{i=1}^s d(P_i) = d(P_\tau) .
$$
Thus we see that $\bx'$ is contained in $\De(P_\tau)$. 
Since $\Int (\tau)$ is a subset of an equivalence class in $\Gamma^*(f)$, 
we have $\De(P_\tau) = \De(P)$ 
for every $P$ in $\Int (\tau)$. 
\end{proof}

\begin{proposition}\label{intersection}
Consider the germ $(f,\zero)$ of a mixed function $f$ and 
let $\Sigma^*$ be a regular simplicial cone subdivision which is admissible for $f$. 
Then, 
for every face $\Delta(P)$, %%
there exist vertices $Q_1, \dots ,Q_s$ in $\Sigma^*$ such that 
$$\Delta(P) \subset \bigcap_{j=1}^s \Delta(Q_j).$$
\end{proposition}

\begin{proof}
Suppose that $\dim \Delta(P) = n-1$. 
There exist vertices $Q_1, \dots ,Q_s$ in $\Sigma^*$ such that 
$\Delta(P) \subset \bigcap_{j=1}^s \Delta(Q_j)$, 
where, of course, we take $s=1$ and $Q_1 = P$. 
Suppose that $\dim \Delta(P) = n-2$. 
Then, in the (radial) dual Newton diagram $\Gamma^*(f)$, 
$P$ is contained in the interior of 
some \underline{admissible (for $\Gamma^*(f)$)} cone $\Xi := \Cone(Q_1, \dots , Q_s)$ 
whose vetices $Q_1, \dots , Q_s$ satisfy $\dim \Delta(Q_j) = n-1$ for every $j = 1, \dots ,s$. 
Hence, $P = \sum_{j=1}^s a_j Q_j$ for some $a_j >0\ (j = 1, \dots ,s)$. %%
By Proposition \ref{non-empty}, we have 
$\De(P) \subset \bigcap_{j=1}^s \De(Q_j)$. 
Here the weight vectors $Q_1, \dots , Q_s$, multiplied by some real positive numbers if necessary, 
are \underline{vertices of $\Sigma^*$}. %%

Suppose that $\dim \Delta(P) = k\ (k \leq n-3)$ and the assertion is true for faces with $\dim \Delta \geq k+1$. 
Then, in the (radial) dual Newton diagram $\Gamma^*(f)$, 
$P$ is contained in the interior of 
some \underline{admissible (for $\Gamma^*(f)$)} cone $\Xi := \Cone(Q_1, \dots , Q_s)$ 
whose vetices $Q_1, \dots , Q_s$ satisfy $\dim \Delta(Q_j) \geq k+1$ for every $j = 1, \dots ,s$. %%
Hence, $P = \sum_{j=1}^s a_j Q_j$ for some $a_j >0\ (j = 1, \dots ,s)$. %%
By Proposition \ref{non-empty}, we have 
$\De(P) \subset \bigcap_{j=1}^s \De(Q_j)$. 
Here, by the induction hypothesis, 
for each face $\De(Q_j)$, 
there exist vertices $Q'_{j1}, \dots ,Q'_{j{s_j}}$ in $\Sigma^*$ such that 
$\Delta(Q_j) \subset \bigcap_{\ell=1}^{s_j} \Delta(Q'_{j \ell})$. 
Thus we have 
$$
\De(P) \subset \bigcap_{j=1}^s \left( \bigcap_{\ell=1}^{s_j} \Delta(Q'_{j \ell}) \right) .
$$
\end{proof}

\medskip

\section{Convenient mixed function germs}\label{section-convenient-function}

Let $(f,\zero)$ be the germ of a mixed function $f$ at $\zero \in \bc^n$. 
For $I \subset \{1, \dots, n\}$, 
$f^I$ denotes the restriction of $f$ 
on the coordinate subspace 
$$\bc^I := \{ \bbz \,|\, z_j=0,\, j\notin I \}.$$

\begin{definition}[\cite{Oka2018}, p.98] \label{convenient-subdivision-f}  %%
{\rm 
Let $\Sigma^*$ be a regular simplicial cone subdivision 
which is admissible for $f$ (Definition \ref{admissible-regul-simpl-cone-subdivision}). 
Then we say $\Sigma^*$ is {\em convenient} 
if for every $I$ with $f^I \not \equiv 0$, the cone $E_{I^c}$ is contained in $\Sigma^*$. 
}
\end{definition}

\begin{definition}[convenient function germ,\ \cite{Oka2018}, p.79]\label{convenient-function} %%
{\rm 
The germ $(f,\zero)$ of a mixed function is called {\em convenient} if 
$f^I \not \equiv 0$ for every $I \neq \emptyset$. %%
}
\end{definition}

\begin{remark}\label{convenient-function-terms}
{\rm 
The germ $(f,\zero)$ of a mixed function $f$ is convenient (Definition \ref{convenient-function}) 
if and only if 
$f$ has some terms 
$c_{\nu, \mu}z_i^{\nu_i}\bar{z}_i^{\mu_i} \ (c_{\nu, \mu}\neq 0)$, where 
$\nu = (0, \dots , 0, \nu_i , 0, \dots , 0),\ \mu = (0, \dots , 0, \mu_i , 0, \dots , 0)$, 
for every $i \ (=1, \dots , n)$. 
}
\end{remark}

\begin{remark} \label{convenient-subdivision-f-important-remark}  %%%
{\rm 
Suppose that 
{\bf (i)}\ $\Sigma^*$ is a convenient (in the sense of Definition \ref{convenient-subdivision-f}) 
regular simplicial cone subdivision which is admissible for $f$ and 
{\bf (ii)}\ the germ of a mixed function $f$ is convenient (Definition \ref{convenient-function}). 
Then $\Sigma^*$ is convenient in the sense of Definition \ref{convenient-subdivision}. 
Hence, 
for the toric modification 
$
\hat{\pi} : X \to \bc^n
$
associated with $\Sigma^*$, 
the restriction 
$
\hat{\pi} : X\setminus \hat{\pi}^{-1}(\zero) \to \bc^n \setminus \{ \zero \} %%%
$
is biholomorphic and we have 
$
\hat{\pi}^{-1}(\zero) = \bigcup_{P\gg 0}\hat{E}(P) %%%
$
 (recall Remark \ref{convenient-subdivision-remark}). 
}
\end{remark}

\begin{theorem}[Isolatedness,\ \cite{Oka2010}, Theorem 19\ (1)]\label{isolated}
Let $(f,\zero)$ be a convenient Newton non-degenerate
\footnote{In Theorem 19 of \cite{Oka2010}, the assertion (1) holds without ``true non-degeneracy" assumption.}
 mixed function germ at $\zero \in \bc^n$. 
Then $\zero$ is 
a mixed regular point of $f$ 
or 
an isolated singular point (Definition \ref{mixed-critical}) of the mixed hypersurface $V:= f^{-1}(0)$. %%
\end{theorem}

\begin{lemma}\label{E-j}
Suppose that $n\geq 2$, the germ $(f,\zero)$ of a mixed function 
$f(\bbz, \bar{\bbz}) = \sum_{\nu, \mu}c_{\nu, \mu} \bbz^\nu \bar{\bbz}^\mu$ 
is convenient (Definition \ref{convenient-function}) and 
$\Sigma^*$ is a convenient (in the sense of Definition \ref{convenient-subdivision-f}) regular simplicial cone subdivision which is admissible for $f$. 
Then we have the following:
\begin{enumerate}
\item All the $1$-dimensional cones $\Cone(P)$ with non-strictly positive vertex $P$ in the (radial) dual Newton diagram $\Gamma^*(f)$ are 
$\Cone(E_1), \dots ,$ and $\Cone(E_n)$. 
\item For a fixed $j \ (= 1, \dots ,n)$, consider the face $\Delta(E_j)   %%
\footnote{$\Delta(E_j)$ is not compact.}.$ 
If $c_{\nu, \mu} \neq 0$ and $\nu + \mu \in \Delta(E_j)$, then we see that 
$
E_j(\nu + \mu) = 0\ \ \text{and}\ \   %%
E_j(\nu - \mu) = 0.$                  %%
\end{enumerate}
\end{lemma}

\begin{proof}
If $n\geq 2$, the germ $(f,\zero)$ of a mixed function $f$ is convenient and 
$\Sigma^*$ is a convenient (in the sense of Definition \ref{convenient-subdivision-f}) regular simplicial cone subdivision which is admissible for $f$, 
then, 
as stated in Remark \ref{convenient-subdivision-f-important-remark}, 
$\Sigma^*$ is convenient in the sense of Definition \ref{convenient-subdivision}. 
Hence, by Remark \ref{convenient-subdivision-remark}, 
all the $1$-dimensional cones $\Cone(P)$ with non-strictly positive vertex $P$ 
in the (radial) dual Newton diagram $\Gamma^*(f)$ are $\Cone(E_1), \dots , \Cone(E_n)$. 
By Remark \ref{convenient-function-terms}, 
$f$ has some terms 
$c_{\nu, \mu}z_i^{\nu_i}\bar{z}_i^{\mu_i} \ (c_{\nu, \mu}\neq 0)$, where 
$\nu = (0, \dots , 0, \nu_i , 0, \dots , 0),\ \mu = (0, \dots , 0, \mu_i , 0, \dots , 0)$, 
for \underline{every} $i \ (=1, \dots , n)$. 
For the weight vector $E_j$ and the term $c_{\nu, \mu}z_i^{\nu_i}\bar{z}_i^{\mu_i}$ with $i\neq j$ and $c_{\nu, \mu}\neq 0$, 
we have 
$$
E_j (\nu + \mu) = \nu_j + \mu_j = 0.  %%
$$
Hence, we see that the minimum value $d(E_j)$ is equal to $0$. 
Thus, we have 
$$
\Delta(E_j) = \{ \xi \in \Gamma_+(f) \ |\ E_j(\xi) = 0 \}    %%
$$
for every $j \ (= 1, \dots ,n)$. 
Thus, 
if $c_{\nu, \mu} \neq 0$ and $\nu + \mu \in \Delta(E_j)$, then we have 
$$E_j(\nu + \mu) = \nu_j + \mu_j = 0,$$    %%
and hence, we have 
$$
\nu_j = \mu_j = 0 .
$$
Thus, moreover, we have
$$
E_j(\nu - \mu) = \nu_j - \mu_j = 0 .   %%
$$
\end{proof}

\medskip

We have an analogy of Proposition 10 (\cite{Oka2015}) for 
mixed function germs of strongly polar \underline{non-negative} mixed weighted homogeneous face type: 
\begin{proposition}\label{change-of-strictly-positive-weight-vector}
Suppose that 
{\bf (1)}\ $n\geq 2$, 
{\bf (2)}\ the germ $(f,\zero)$ of a mixed function 
$f(\bbz, \bar{\bbz}) = \sum_{\nu, \mu}c_{\nu, \mu} \bbz^\nu \bar{\bbz}^\mu$ 
is convenient, and 
{\bf (3)}\ $\Sigma^*$ is a convenient (Definition \ref{convenient-subdivision-f}) regular simplicial cone subdivision which is admissible for $f$. 
If the germ $(f,\zero)$ of a mixed function $f$ 
is of strongly polar non-negative mixed weighted homogeneous face type (Definition \ref{strongly-mixed-wh-face-type}), 
then 
for any strictly positive weight vector $P$, 
$f_P$ is a strongly polar non-negative mixed weighted homogeneous polynomial
with respect to the weight vector $P$. 
\end{proposition}

\begin{proof}
By the definition, for every compact face $\Delta$, 
the face function $f_\Delta(\bbz, \bar{\bbz})$ is 
a strongly polar non-negative mixed weighted homogeneous polynomial (Definition \ref{strongly-mixed-whp}). 

Suppose that $\dim \Delta(P) = n-1$. By the assumption, %%%
$f_P$ is a strongly polar non-negative mixed weighted homogeneous polynomial with respect to $P$ 
since $P$ is unique up to multiplications of positive real numbers. %%%
Suppose that $\dim \Delta(P) = n-2$. 
Then, in the (radial) dual Newton diagram $\Gamma^*(f)$, 
$P$ is contained in the interior of 
some admissible (for $\Gamma^*(f)$) cone $\Xi := \Cone(Q_1, \dots , Q_s)$ 
whose vetices $Q_1, \dots , Q_s$ satisfy $\dim \Delta(Q_j) = n-1$ for every $j = 1, \dots ,s$. 
We have $P = \sum_{j=1}^s a_j Q_j$ for some $a_j >0\ (j = 1, \dots ,s)$. 
By Proposition \ref{non-empty}, we have 
$$
\De(P) \subset \bigcap_{j=1}^s \De(Q_j).
$$
Here the weight vectors $Q_1, \dots , Q_s$, multiplied by some real positive numbers if necessary, 
are vertices of $\Sigma^*$. 
If $c_{\nu, \mu} \neq 0$ and $\nu + \mu \in \De(P)$, 
then we have 
$\nu + \mu \in \De(Q_j)$ for every $j = 1, \dots ,s$. 
By the assumptions of this Proposition, 
all the non-strictly positive vertices of $\Sigma^*$ are $E_1, \dots , E_n$. 
If $c_{\nu, \mu} \neq 0$ and $\nu + \mu \in \Delta(E_j)$, 
then, by Lemma \ref{E-j}, we have 
$
E_j(\nu + \mu) = 0\ \ \text{and}\ \  %%
E_j(\nu - \mu) = 0                   %%
$. 
Hence, we may say that 
for every $(\nu, \mu)$ with $c_{\nu, \mu} \neq 0$ and $\nu + \mu \in \De(P)$, %%
$$
P(\nu - \mu) = \sum_{j=1}^s a_j Q_j(\nu - \mu) = \sum_{j=1}^t a_j Q_j(\nu - \mu) , %%
$$
where $Q_1, \dots ,Q_t$ are all strictly positive. 
By the above result on $(n-1)$-dimensional case, 
the numbers $Q_1(\nu - \mu), \dots ,Q_t(\nu - \mu)$ are all non-negative, %%
and each number $m_j := Q_j(\nu - \mu)$ does not depend on $(\nu, \mu)$ with $c_{\nu, \mu} \neq 0$ and $\nu + \mu \in \De(P)$. 
Thus we conclude that 
$$\pdeg_P f_P = P(\nu - \mu) = \sum_{j=1}^t a_j m_j \geq 0. $$

Now suppose that $\dim \Delta(P) = k\ (k \leq n-3)$ and the assertion is true for faces with $\dim \Delta \geq k+1$. 
In the (radial) dual Newton diagram $\Gamma^*(f)$, 
$P$ is contained in the interior of 
some admissible (for $\Gamma^*(f)$) cone $\Xi := \Cone(Q_1, \dots , Q_s)$ 
whose vetices $Q_1, \dots , Q_s$ satisfy $\dim \Delta(Q_j) \geq k+1$ for every $j = 1, \dots ,s$. 
We have $P = \sum_{j=1}^s a_j Q_j$ for some $a_j >0\ (j = 1, \dots ,s)$. 
By Proposition \ref{non-empty}, we have 
$\De(P) \subset \bigcap_{j=1}^s \De(Q_j)$. 
Here, by using Proposition \ref{intersection}, 
we may say that all $Q_1, \dots , Q_s$ are vetices of $\Sigma^*$. 
By the assumptions of this Proposition, 
all the non-strictly positive vertices of $\Sigma^*$ are $E_1, \dots , E_n$. 
If $c_{\nu, \mu} \neq 0$ and $\nu + \mu \in \Delta(E_j)$, 
then, by Lemma \ref{E-j}, we have 
$
E_j(\nu + \mu) = 0\ \ \text{and}\ \  %%
E_j(\nu - \mu) = 0                   %%
$. 
Hence, we may say that for every $(\nu, \mu)$ with $c_{\nu, \mu} \neq 0$ and $\nu + \mu \in \De(P)$, 
$$
P(\nu - \mu) = \sum_{j=1}^s a_j Q_j(\nu - \mu) = \sum_{j=1}^t a_j Q_j(\nu - \mu) , %%
$$
where $Q_1, \dots ,Q_t$ are all strictly positive. 
By the induction hypothesis, 
the numbers $Q_1(\nu - \mu), \dots ,Q_t(\nu - \mu)$ are all non-negative, 
and each number $m_j := Q_j(\nu - \mu)$ does not depend on $(\nu, \mu)$ with $c_{\nu, \mu} \neq 0$ and $\nu + \mu \in \De(P)$. 
Thus we see that 
$$\pdeg_P f_P = P(\nu - \mu) = \sum_{j=1}^t a_j m_j \geq 0. $$
\end{proof}

\medskip

\section{Resolutions of Newton non-degenerate mixed polynomials of 
strongly polar non-negative mixed weighted homogeneous face type}\label{section-main}

Theorem 11 (\cite{Oka2015}) and Theorem 9.19 (\cite{Oka2018}) 
are generalizations of Varchenko's famous Lemma 10.3 (\cite{Varchenko1976}) on holomorphic function germs 
to mixed polynomial germs or mixed function germs. 

\bigskip

Our main Theorem \ref{theorem11-improved} below improves 
Theorem 11 (1), (2) of \cite{Oka2015} under the {\sc Assumption}{\rm (*)}. 
In order to prove Theorem \ref{theorem11-improved}, 
we prepare the following Lemma \ref{fraction_lemma}. 

\medskip

Here we say a complex valued function $\rho(\bbu)$ with complex $n$ variables 
$\bbu = (u_1, \dots ,u_n), u_j = x_j + iy_j \ (j = 1, \dots ,n)$ 
is {\em of class $C^s$} 
if the real and imaginary parts of 
$\rho(\bbu)= \rho(x_1,y_1, \dots ,x_n,y_n)$ are $C^s$-functions with real $2n$ variables. 

If the real and imaginary parts of $\rho(\bbu)= \rho(x_1,y_1, \dots ,x_n,y_n)$ are of class $C^1$, 
namely, they have their partial derivatives 
$\frac{\partial \Re \rho}{\partial x_j},\ \frac{\partial \Re \rho}{\partial y_j}$, 
$\frac{\partial \Im \rho}{\partial x_j},\ \frac{\partial \Im \rho}{\partial y_j}$ 
and these $4$ functions are continuous, 
then by \eqref{partial1}, 
$\frac{\partial \rho}{\partial z_j}$ and $\frac{\partial \rho}{\partial \bar{z}_j}$ exist 
and they are continuous. 
Conversely, if $\frac{\partial \rho}{\partial z_j}$ and $\frac{\partial \rho}{\partial \bar{z}_j}$ exist 
and they are continuous, 
then by \eqref{partial2}, 
the real and imaginary parts of $\rho(\bbu)= \rho(x_1,y_1, \dots ,x_n,y_n)$ are of class $C^1$. 

\begin{lemma}\label{fraction_lemma}
{\rm 
Let $r,\ s\ (\geq 1)$ be integers, and $u$ be a complex variable. We set
$$
\xi(u) := 
\begin{cases}
\displaystyle \frac{u^{r+s}}{\bar{u}^r} & (u \neq 0)\\
0                         & (u = 0) .
\end{cases}
$$
Then $\xi(u)$ is of class $C^{s-1}$ on $\bc$. Hence, 
$$
\bar{\xi}(u) = 
\begin{cases}
\displaystyle \frac{\bar{u}^{r+s}}{u^r} & (u \neq 0)\\
0                         & (u = 0)
\end{cases}
$$
is also of class $C^{s-1}$ on $\bc$. 
}
\end{lemma}

\begin{proof}
We have
\begin{equation}\label{re_and_im}
\xi(u) = \xi(x,y) = 
\displaystyle 
\frac{u^{r+s}}{\bar{u}^r} = \frac{u^{{r+s}+r}}{|u|^{2r}}
= \frac{(x+iy)^{2r+s}}{\sqrt{x^2+y^2 \,}^{2r}} = \frac{(x+iy)^{2r+s}}{(x^2+y^2)^r}\ \ \ (\text{if}\ (x,y)\neq (0,0)) .
\end{equation}
The real and imaginary parts of the last numerator of \eqref{re_and_im} are of degree $2r+s$. 
We will prove that the real and imaginary parts of $\xi(u)=\xi(x,y)$ are of class $C^{s-1}$. 
To do this, it is sufficient that we show 
$$
\varphi(x,y):= 
\begin{cases}
\displaystyle \frac{x^{2r+s-j}y^j}{(x^2+y^2)^r} & ((x,y) \neq (0,0))\\
0                         & ((x,y) = (0,0))
\end{cases}
$$
is a $C^{s-1}$-function on $\br^2$ for every $r\ (\geq 1)$ and every $j\ (0\leq j\leq 2r+s)$. 
Therefore we will prove this assertion by the induction on $s$. 

If $s=1$, then 
$$
\displaystyle |\varphi(x,y)| = \frac{|x|^{2r+1-j}|y|^j}{\sqrt{x^2+y^2 \,}^{2r}}
\leq |x|\ \text{or}\ |y|
\ \ \ \ (\text{if}\ (x,y)\neq (0,0)) .
$$
Hence, $\varphi(x,y)$ is a continuous, namely, $C^0$-function on $\br^2$. 
Thus the assertion is true for $s=1$. 

Suppose that the assertion is true 
for $s$ (with $s\geq 1$), every $r\ (\geq 1)$ and every $j\ (0\leq j\leq 2r+s)$. 
We will prove that 
$$
\psi(x,y):= 
\begin{cases}
\displaystyle \frac{x^{2r+s+1-j}y^j}{(x^2+y^2)^r} & ((x,y) \neq (0,0))\\
0                         & ((x,y) = (0,0))
\end{cases}
$$
is a $C^s$-function on $\br^2$ for every $r\ (\geq 1)$ and $j\ (0\leq j\leq 2r+s+1)$. 

We have
$$
\begin{array}{ccl}
\displaystyle \frac{\partial \psi}{\partial x}(x,y)& = & 
\frac{(2r+s+1-j)x^{2r+s-j}y^j(x^2+y^2)^r - 2rx\cdot x^{2r+s+1-j}y^j(x^2+y^2)^{r-1}}{(x^2+y^2)^{2r}}\\
  & = & \frac{(2r+s+1-j)x^{2r+s-j}y^j}{(x^2+y^2)^r} - 2r\frac{x^{2(r+1)+s-j}y^j}{(x^2+y^2)^{r+1}}\ \ \ ((x,y) \neq (0,0))
\end{array}
$$
and 
$$
\begin{array}{ccl}
\displaystyle \frac{\partial \psi}{\partial y}(x,y)& = & 
\frac{jx^{2r+s+1-j}y^{j-1}(x^2+y^2)^r - 2ry\cdot x^{2r+s+1-j}y^j(x^2+y^2)^{r-1}}{(x^2+y^2)^{2r}}\\
  & = & \frac{jx^{2r+s-(j-1)}y^{j-1}}{(x^2+y^2)^r} - 2r\frac{x^{2(r+1)+s-(j+1)}y^{j+1}}{(x^2+y^2)^{r+1}}\ \ \ ((x,y) \neq (0,0)) .
\end{array}
$$

Since 
$$
\psi(x,0)= x^{s+1-j}y^j = 
\begin{cases}
x^{s+1} & (j=0) \\
0       & (j\geq 1)
\end{cases}
$$
and 
$$
\psi(0,y)= x^{2r+s+1-j}y^{j-2r} = 
\begin{cases}
0       & (j\leq 2r+s)\\
y^{s+1} & (j = 2r+s+1) ,
\end{cases}
$$
we have 
$$
\displaystyle \frac{\partial \psi}{\partial x}(0,0) = \frac{\partial \psi}{\partial y}(0,0) = 0 .
$$
Thus, by the induction hypothesis, 
$\frac{\partial \psi}{\partial x}(x,y)$ and $\frac{\partial \psi}{\partial y}(x,y)$ 
are $C^{s-1}$-functions on $\br^2$. 
Hence, $\psi(x,y)$ is a $C^s$-function on $\br^2$. 
This completes the proof of Lemma \ref{fraction_lemma}. 
%%%%%%%%%%%%%%%%%%%%%%%%%%%%%%%%%%%%%%%%%%%%%%%%%%%%%

\medskip

We can give another proof using Wirtinger derivatives $\frac{\partial}{\partial u},\ \frac{\partial}{\partial \bar{u}}$. 
We prove the assertion of Lemma \ref{fraction_lemma} by the induction on $s$. 
If $s=1$, then $\xi(u)$ is a real analytic map on $\bc \setminus \{ 0 \}$. 
Since $|\xi(u)|= |u|\ (\text{if}\ u\neq 0)$, we have $\lim_{u \to 0}\xi(u) = 0 = \xi(0)$. 
Hence, $\xi(u)$ is a continuous, namely, $C^0$-map on $\bc$. 
Thus the assertion is true for $s=1$. 
Suppose that the assertion is true for $s$ (with $s\geq 1$) and every $r\ (\geq 1)$. 
We will prove 
$$
\eta(u) := 
\begin{cases}
\displaystyle \frac{u^{r+s+1}}{\bar{u}^r} & (u \neq 0)\\
              0                           & (u = 0)
\end{cases}
$$
is of class $C^s$ on $\bc$ for every $r\ (\geq 1)$.
We have
\begin{equation}\label{partial_holo}
\frac{\partial \eta}{\partial u}(u) = \frac{(r+s+1)u^{r+s}\bar{u}^r - u^{r+s+1} \cdot 0}{\bar{u}^{2r}}
= (r+s+1)\frac{u^{r+s}}{\bar{u}^r}\ \ (\text{if}\ u\neq 0)
\end{equation}
and
\begin{equation}\label{partial_conj}
\frac{\partial \eta}{\partial \bar{u}}(u) = \frac{0 \cdot \bar{u}^r - ru^{r+s+1}\bar{u}^{r-1}}{\bar{u}^{2r}}
= -r\frac{u^{r+1+s}}{\bar{u}^{r+1}}\ \ (\text{if}\ u\neq 0) .
\end{equation}
We set $u = x + iy = (x,y)\ (x,y\in \br)$. 
Since $\eta(x,0)= x^{s+1}$ and $\eta(0,y)= (iy)^{r+s+1}(-iy)^{-r}= (-1)^ri^{s+1}y^{s+1}$, we have
$$
\frac{\partial \eta}{\partial x}(0,0) = 0\ \ \text{and}\ \ 
\frac{\partial \eta}{\partial y}(0,0) = 0. 
$$
Hence, we have 
\begin{equation}\label{partial_holo_zero}
\frac{\partial \eta}{\partial u}(0) = 0
\end{equation}
and
\begin{equation}\label{partial_conj_zero}
\frac{\partial \eta}{\partial \bar{u}}(0) = 0 .
\end{equation}

From \eqref{partial_holo} and \eqref{partial_holo_zero}, 
by the induction hypothesis, 
we see $\frac{\partial \eta}{\partial u}$ is of class $C^{s-1}$ on $\bc$. 
From \eqref{partial_conj} and \eqref{partial_conj_zero}, 
%% by the induction hypothesis, 省略
we also see $\frac{\partial \eta}{\partial \bar{u}}(u)$ is of class $C^{s-1}$ on $\bc$. 
Hence, by \eqref{partial2}, $\frac{\partial \eta}{\partial x}$ and $\frac{\partial \eta}{\partial y}$ 
are also of class $C^{s-1}$ on $\bc$. 
Thus it is proved that $\eta$ is of class $C^s$ on $\bc$. 
This completes the proof of Lemma \ref{fraction_lemma}. 
%%%%%%%%%%%%%%%%%%%%%%%%%%%%%%%%%%%%%%%%%%%%%%%%%%%%%%%%
\end{proof}

\vspace{1cm}

We now set up the situation of Theorem \ref{theorem11-improved}. 
Let $(f,\zero)$ be 
a convenient          
Newton non-degenerate 
mixed polynomial germ 
$$
f(\bbz, \bar{\bbz}) = \sum_{\nu, \mu}c_{\nu, \mu} \bbz^\nu \bar{\bbz}^\mu 
$$
at $\zero \in \bc^n$ 
of strongly polar non-negative mixed weighted homogeneous face type (Definition \ref{strongly-mixed-wh-face-type}). 
We assume $f(\zero)=0$. 
Then by Theorem \ref{isolated}, 
$\zero$ is a mixed regular point of $f$ 
or 
an isolated singular point (Definition \ref{mixed-critical}) of the mixed hypersurface $V:= f^{-1}(0)$. 

\medskip

Let $\Sigma^*$ be a regular simplicial subdivision (of $N^+_\br$) which is admissible for $f$ and 
assume that $\Sigma^*$ is convenient in the sense of Definition \ref{convenient-subdivision-f} 
(recall Remark \ref{convenient-subdivision-f-important-remark}). 
Consider the toric modification $\hat{\pi} : X \to \bc^n$ associated with $\Sigma^*$. 
Then (recall Remark \ref{convenient-subdivision-remark}) 
the toric modification 
$$
\hat{\pi} : X \to \bc^n
$$
associated with $\Sigma^*$ 
is a proper birational map and 
$
\hat{\pi} : X\setminus \hat{\pi}^{-1}(\zero) \to \bc^n \setminus \{ \zero \}
$
is biholomorphic and we have 
$$
\hat{\pi}^{-1}(\zero) = \bigcup_{P\gg 0}\hat{E}(P),
$$
where $\hat{E}(P)$ is the exceptional divisor defined by a vertex $P$ of $\Sigma^*$. 

Let $(V,\zero)$ be the germ of the mixed hypersurface $V := f^{-1}(0)$ at $\zero \in \bc^n$. 
We call 
$$\tilde V := \text{the closure of}\ \hat{\pi}\inv(V\setminus\{ \zero \})\ \text{in}\ X$$
the {\em strict transform of $V$ to $X$}. 
We consider the strict transform $\tilde V$ 
as the germ at $\tilde V \cap \hat{\pi}^{-1}(\zero)$. 
For each $k$-dimensional cone $\tau = \Cone(P_1, \dots , P_k)$ of $\Sigma^*$, 
we set 
$$
\hat{E}(\tau)^* := (\cap_{i=1}^k \hat{E}(P_i)) \setminus (\cup_{Q\in \Va, Q\not\in \tau}\hat{E}(Q))
$$
and 
$$
\tilde V(\tau)^* := \tilde V \cap \hat{E}(\tau)^* . 
$$

\medskip

If the vertices of $\tau = \Cone(P_1, \dots , P_k)$ of $\Sigma^*$ 
are all \underline{strictly positive} 
and 
$$\{ ((\nu, \mu),\ j)\ |\ c_{\nu, \mu}\neq 0,\ \nu + \mu \not\in \De(P_j),\ j = 1, \dots, k \} \neq \emptyset,$$
where $r_j := \rdeg_{P_j} f_{P_j} = d(P_j)$, 
then we define 
$$
\Lambda(\tau):= \min \{ P_j(\nu + \mu) - r_j\ |\ c_{\nu, \mu}\neq 0,\ \nu + \mu \not\in \De(P_j),\ j = 1, \dots, k \}. 
$$
By the definition, we have $\Lambda(\tau) \geq 1$. 
Moreover, we define $L(\Sigma^*)$ 
to be \underline{the set of all $\Lambda(\tau)$'s} such that \\
\ \ \ \ {\rm (1)}\ $\tau$ is a cone of $\Sigma^*$ whose vertices are all strictly positive,\\
\ \ \ \ {\rm (2)}\ $\{ ((\nu, \mu),\ j)\ |\ c_{\nu, \mu}\neq 0,\ \nu + \mu \not\in \De(P_j),\ j = 1, \dots, k \} \neq \emptyset$, where $k:= \dim \tau$, and \\
\ \ \ \ {\rm (3)}\ $\tilde V(\tau)^* \cap \bc_\sigma^n \neq \emptyset$ for some $n$-dimensional cone $\sigma$ of $\Sigma^*$ with $\tau \prec \sigma$. 

\medskip

\noindent
If $L(\Sigma^*) \neq \emptyset$, then we define 
$$
\La := \min L(\Sigma^*) .
$$
By the definition, we have $\La \geq 1$. 

\medskip

\begin{theorem}\label{theorem11-improved}
We assume $n\geq 2$. 
Under the situation mentioned above, we have the following:

\medskip

\noindent
{\bf (i)} If the vertices of the $k$-dimensional cone $\tau$ of $\Sigma^*$ are all strictly positive, then 
$\tilde V(\tau)^*$ is empty or a real analytic manifold of dimension $2(n-k-1)$. 

\medskip

\noindent
{\bf (ii)} We assume:

\medskip

\ \ {\sc Assumption}{\rm (*)}:\ 
For every point $\hat{\bbu}^0 \in \tilde V \cap \hat{\pi}^{-1}(\zero) = \tilde V \cap \bigcup_{P\gg 0}\hat{E}(P)$, 
there exists an $n$-dimensional 

\ cone $\sigma = (P_1, \dots , P_n)$ 
with strictly positive vertices $P_1, \dots , P_\kappa \ (1 \leq \kappa \leq n)$ 
for which the toric chart 

\ $\bc_\sigma^n$ contains $\hat{\bbu}^0$ 
and 
$\hat{\bbu}^0 = \bbu_\si^0 = 
(\underbrace{0, \dots, 0}_{\kappa} \, ,\ u^0_{\kappa +1}, \dots, u^0_n),\ \ u^0_{\kappa +1}\neq 0, \dots ,u^0_n\neq 0
\ \ \ \text{in}\ \bc_\sigma^n$. 

\medskip

\noindent
Then the strict transform germ $\tilde V$ at $\tilde V \cap \hat{\pi}^{-1}(\zero)$ 
is 
a $C^{\La - 1}$-manifold if $L(\Sigma^*) \neq \emptyset$, and 
a real analytic manifold if $L(\Sigma^*) = \emptyset$. 
Hence, it is a topological manifold. 
Moreover, it is a real analytic manifold outside of $\tilde V \cap \hat{\pi}^{-1}(\zero)$. 
\end{theorem}

\medskip

\begin{proof}
We first prove the assertion {\bf (i)}. 
Let $\tau := (P_1, \dots , P_k)$ be a $k$-dimensional cone in $\Sigma^*$ whose vertices are all strictly positive vectors. 
Then the faces $\De(P_1), \dots , \De(P_k)$ of $\Gamma_+(f)$ are all compact. 
If $\tilde V(\tau)^* = \emptyset$, then the proof is completed, and 
if $\tilde V(\tau)^* \neq \emptyset$, then we can take an $n$-dimensional cone
$$\sigma = (P_1, \dots , P_n)$$
of $\Sigma^*$ and let 
$$\bbu_\si = (u_{\si 1},\dots, u_{\si n})$$ 
be the toric coordinates of the chart $\bc_\sigma^n$. 
We write $u_{\si j}=u_j$ hereafter. 

We put 
$$
\De := \bigcap_{j=1}^k \Delta(P_j) .
$$
Then $\De$ is a compact face of $\Gamma_+(f)$ (recall Proposition \ref{intersection-of-faces}
\footnote{$\De = \De(P_\tau)$ where we set $P_\tau := P_1 + \cdots + P_k$. Hence, $P_\tau$ is a strictly positive vector.}
). 
For each $j=1,\dots, k$, we set 
$$r_j := \rdeg_{P_j} f_{P_j} , \ \ \ \text{and} \ \ \ p_j := \pdeg_{P_j } f_{P_j}.$$
By the strongly polar non-negative mixed weightedness assumption 
and Proposition \ref{change-of-strictly-positive-weight-vector}, 
we have $r_j > 0$ and $p_j \geq 0$ for every  $j \ (= 1, \dots , k)$.

Put 
$$P_j = {}^t(p_{1j}, \dots , p_{nj}).$$
If $\nu + \mu \in \Delta(P_j)$, then we have 
$$
\sum_{i=1}^n p_{ij}(\nu_j + \mu_j) = r_j >0\ \ \ \ \text{and}\ \ \ \ 
\sum_{i=1}^n p_{ij}(\nu_i - \mu_i) =  p_j \geq 0 .
$$
Hence, note that $r_j + p_j >0$ and %%
$r_j + p_j = 2\sum_{i=1}^n p_{ij}\nu_j$ is an even positive integer. On the other hand, we have 
$r_j - p_j = 2\sum_{i=1}^n p_{ij}\mu_j$ is an even non-negative integer. 
Hence, we have two integers
$$
\frac{r_j + p_j}{2} >0,\ \ \ \ \frac{r_j - p_j}{2} \geq 0 .
$$

Recall that $\hat{\pi}$ coincides with the toric morphism $\pi_\si$ on the toric chart $\bc_\sigma^n$. 
Let us consider the pull-back polynomial 
$$
\pi_\si^*f := f \circ \pi_\si . 
$$
We decompose the polynomial $\pi_\si^*f$ ``formally" as below: %%
$$
\begin{array}{cl}
  & \pi_\si^*f(\bbu,\bar \bbu) \\
= & \displaystyle \sum_{\nu, \mu}c_{\nu, \mu} u_1^{P_1(\nu)} \cdots u_n^{P_n(\nu)} \bar{u}_1^{P_1(\mu)} \cdots \bar{u}_n^{P_n(\mu)} \\
= & \displaystyle \prod_{j=1}^k u_j^{\frac{r_j+p_j}2}\bar u_ j^{\frac{r_j-p_j}2} \, \times 
\left( \sum_{\nu, \mu}c_{\nu, \mu} 
\prod_{j=1}^k u_j^{P_j(\nu) - \frac{r_j+p_j}2} \prod_{j=1}^k \bar{u}_j^{P_j(\mu) - \frac{r_j-p_j}2} 
u_{k+1}^{P_{k+1}(\nu)} \cdots u_n^{P_n(\nu)} \bar{u}_{k+1}^{P_{k+1}(\mu)} \cdots \bar{u}_n^{P_n(\mu)} \right),
\end{array}
$$
and we set 
$$
\widetilde f(\bbu,\bar \bbu) := \displaystyle \sum_{\nu, \mu}c_{\nu, \mu} 
\prod_{j=1}^k u_j^{P_j(\nu) - \frac{r_j+p_j}2} \prod_{j=1}^k \bar{u}_j^{P_j(\mu) - \frac{r_j-p_j}2} 
u_{k+1}^{P_{k+1}(\nu)} \cdots u_n^{P_n(\nu)} \bar{u}_{k+1}^{P_{k+1}(\mu)} \cdots \bar{u}_n^{P_n(\mu)} .
$$
For every $j \ (= 1, \dots ,k)$, we have 
$$
u_j^{P_j(\nu)} \bar{u}_j^{P_j(\mu)} = u_j^{P_j(\nu + \mu)/2 + P_j(\nu - \mu)/2} \bar{u}_j^{P_j(\nu + \mu)/2 - P_j(\nu - \mu)/2} .
$$
If $\nu + \mu \in \De$, then we have 
$$
u_j^{P_j(\nu)} \bar{u}_j^{P_j(\mu)} = u_j^{\frac{r_j+p_j}2} \bar u_ j^{\frac{r_j-p_j}2}  %%
$$
for every $j \ (= 1, \dots ,k)$ 
by the strongly mixed weightedness assumption. %%
Put
$$
\begin{array}{ccl}
\widetilde f_\De & := &   %%%%
\displaystyle \sum_{\nu + \mu \in \De}c_{\nu, \mu} 
\prod_{j=1}^k u_j^{P_j(\nu) - \frac{r_j+p_j}2} \prod_{j=1}^k \bar{u}_j^{P_j(\mu) - \frac{r_j-p_j}2} 
u_{k+1}^{P_{k+1}(\nu)} \cdots u_n^{P_n(\nu)} \bar{u}_{k+1}^{P_{k+1}(\mu)} \cdots \bar{u}_n^{P_n(\mu)} \\
 & = & 
\displaystyle \sum_{\nu + \mu \in \De}c_{\nu, \mu} \,
u_{k+1}^{P_{k+1}(\nu)} \cdots u_n^{P_n(\nu)} \bar{u}_{k+1}^{P_{k+1}(\mu)} \cdots \bar{u}_n^{P_n(\mu)} .
\end{array}
$$
Note that 
$\widetilde f_\De$ is a mixed \underline{polynomial} 
having only the variables $u_{k+1},\dots,u_n$, 
namely, 
$\widetilde f_\De$ does not contain any variables $u_{ 1},\dots, u_{ k}$). %%

Thus we obtain the following ``formal" expression: %%
\begin{equation}\label{reminder}
\widetilde f(\bbu,\bar \bbu) = \widetilde f_\De(\bbu',\bar\bbu') + \widetilde R(\bbu,\bar\bbu)
\end{equation}
where $\bbu':=(u_{k+1},\dots,u_{n})$ and 
$$
\widetilde R(\bbu,\bar\bbu) := \displaystyle \sum_{\nu + \mu \not\in  \De} 
c_{\nu, \mu} \prod_{j=1}^k u_j^{P_j(\nu) - \frac{r_j+p_j}2} \prod_{j=1}^k \bar{u}_j^{P_j(\mu) - \frac{r_j-p_j}2} 
u_{k+1}^{P_{k+1}(\nu)} \cdots u_n^{P_n(\nu)} \bar{u}_{k+1}^{P_{k+1}(\mu)} \cdots \bar{u}_n^{P_n(\mu)} .
$$
If there are no $(\nu, \mu)$ with $c_{\nu, \mu} \neq 0$ and $\nu + \mu \not\in \De$, 
then $\widetilde R(\bbu,\bar\bbu) \equiv 0$. 

Note that it is possible that 
some of $P_j(\nu) - \frac{r_j+p_j}2,\ P_j(\mu) - \frac{r_j-p_j}2\ (j = 1, \dots, k)$ are \underline{negative} integers. %%
On the other hand, 
$P_j(\nu) \geq 0$ and $P_j(\mu) \geq 0$ for every $j \ (=k+1, \dots , n)$. 
Set 
$$
\begin{array}{lll}
a_j(\nu) := P_j(\nu) - \frac{r_j+p_j}2, \ \ & b_j(\mu) := P_j(\mu) - \frac{r_j-p_j}2\ \ & \text{for each}\ j \ (= 1, \dots, k),\ \text{and}\\
a_j(\nu) := P_j(\nu),\ \                       & b_j(\mu) := P_j(\mu)\ \                      & \text{for each}\ j \ (= k+1, \dots, n). 
\end{array}
$$
Then we have another expression: 
$$
\widetilde R(\bbu,\bar\bbu) = \displaystyle \sum_{\nu + \mu \not\in  \De} 
c_{\nu, \mu} \left( \prod_{j=1}^k u_j^{a_j(\nu)}\bar{u}_j^{b_j(\mu)} \right) \cdot 
u_{k+1}^{a_{k+1}(\nu)} \cdots u_n^{a_n(\nu)} \bar{u}_{k+1}^{b_{k+1}(\mu)} \cdots \bar{u}_n^{b_n(\mu)} .
$$

\bigskip

We now discuss the continuity of $\widetilde R(\bbu,\bar\bbu)$. 
Let us pay attention to each ``monomial"
$$
B_{\nu, \mu}(\bbu,\bar\bbu) := 
c_{\nu, \mu} \left( \prod_{j=1}^k u_j^{P_j(\nu) - \frac{r_j+p_j}2} \bar{u}_j^{P_j(\mu) - \frac{r_j-p_j}2} \right) \cdot 
u_{k+1}^{P_{k+1}(\nu)} \cdots u_n^{P_n(\nu)} \bar{u}_{k+1}^{P_{k+1}(\mu)} \cdots \bar{u}_n^{P_n(\mu)},\ \ 
c_{\nu, \mu} \neq 0
$$
of $\widetilde R(\bbu,\bar\bbu)$. 

For each $j \ (= 1, \dots, k)$, we have 
$$
a_j(\nu) + b_j(\mu) = P_j(\nu + \mu) - r_j .
$$

If $(\nu, \mu)$ satisfies $\nu + \mu \not\in  \De$, 
then there exists $j \ (= 1, \dots, k)$ such that $\nu + \mu \not\in  \De(P_j)$, 
and for such $j$, we have 
$$
a_j(\nu) + b_j(\mu) = P_j(\nu + \mu) - r_j \geq 1 . %%
$$
Then, by Lemma \ref{fraction_lemma}, 
the function $u_j^{a_j(\nu)}\bar{u}_j^{b_j(\mu)}$ is of class $C^{(P_j(\nu + \mu) - r_j) - 1}$ on $\bc_\sigma^n$ 
if we define the value of the function $u_j^{a_j(\nu)}\bar{u}_j^{b_j(\mu)}$ at $u_j = 0$ to be $0$. 
Especially, 
we see that $u_j^{a_j(\nu)}\bar{u}_j^{b_j(\mu)}$ is a continuous function on $\bc_\sigma^n$. 

On the other hand, if $(\nu, \mu)$ satisfies $\nu + \mu \in \De(P_j)$ for $j \ (= 1, \dots, k)$,
then we have 
$a_j(\nu) + b_j(\mu) = P_j(\nu + \mu) - r_j =0$, namely, 
$P_j(\nu + \mu) = r_j$, 
and also have 
$P_j(\nu - \mu) = p_j$ 
since the face function $f_{P_j}$ is a polar weighted homogeneous polynomial of polar degree $p_j$ with respect to $P_j$. 
Thus, for such $j$, we have 
$$
a_j(\nu) = b_j(\mu) = 0 ,
$$
and hence, corresponding ``monomial" 
$B_{\nu, \mu}(\bbu,\bar\bbu)$ of $\widetilde R(\bbu,\bar\bbu)$ 
does not contain such variables $u_j$ and $\bar{u}_j$. 

Any way, for each ``monomial" $B_{\nu, \mu}(\bbu,\bar\bbu)$ of $\widetilde R(\bbu,\bar\bbu)$, 
there \underline{exists} $j \ (= 1, \dots, k)$ such that $\nu + \mu \not\in  \De(P_j)$ as stated above, %%
and hence, we finally define 
$$
B_{\nu, \mu}(\bbu,\bar\bbu) = 0
$$
when 
$u_j = 0$ where $j = 1, \dots , k$ and $\nu + \mu \not\in  \De(P_j)$.  %%
For fixed $(\nu, \mu)$ with $c_{\nu, \mu} \neq 0$ and $\nu + \mu \not\in \De$, %%
we set 
$$
\lambda_{\nu, \mu}(\tau) := \min \{ P_j(\nu + \mu) - r_j \ |\  \nu + \mu \not\in  \De(P_j),\ j = 1, \dots, k \} ,
$$
where we know that $\{ P_j(\nu + \mu) - r_j \ |\  \nu + \mu \not\in  \De(P_j),\ j = 1, \dots, k \} \neq \emptyset$. 
Since $P_j(\nu + \mu) - r_j \geq 1$ for every $j \ (= 1, \dots ,k)$ with $\nu + \mu \not\in  \De(P_j)$, 
we see that $\lambda_{\nu, \mu}(\tau) \geq 1$. 
Thus, it is concluded that 
each ``monomial" $B_{\nu, \mu}(\bbu,\bar\bbu)$ of $\widetilde R(\bbu,\bar\bbu)$ is of class $C^{\lambda_{\nu, \mu}(\tau) -1}$, 
and hence, continuous on $\bc_\si^n$. 

\medskip

Recall that 
if there are no $(\nu, \mu)$ with $c_{\nu, \mu} \neq 0$ and $\nu + \mu \not\in \De$, 
then $\widetilde R(\bbu,\bar\bbu) \equiv 0$. 
If there exists $(\nu, \mu)$ with $c_{\nu, \mu} \neq 0$ and $\nu + \mu \not\in \De$, 
then we see $\{ \lambda_{\nu, \mu}(\tau)\ |\ c_{\nu, \mu}\neq 0,\ \nu + \mu \not\in \De \} \neq \emptyset$, 
and we conclude that
$$
\Lambda(\tau) = \min \{ \lambda_{\nu, \mu}(\tau)\ |\ c_{\nu, \mu}\neq 0,\ \nu + \mu \not\in \De \} .
$$
Note that $\Lambda(\tau) \geq 1$. 
Thus, 
$
\widetilde R(\bbu,\bar\bbu) = \displaystyle \sum_{\nu + \mu \not\in  \De} B_{\nu, \mu}(\bbu,\bar\bbu)
$
is of class $C^{\Lambda(\tau) - 1}$, and hence, continuous on $\bc_\si^n$. 
Note that 
\begin{equation}\label{widetildeR0dots0}
\widetilde R(\underbrace{0, \dots , 0}_{k} \, , \bbu', \bar\bbu') = 0 .
\end{equation}

We finally define the function 
$\widetilde f(\bbu,\bar \bbu)$ 
by the equality \eqref{reminder} again. 
Then, $\widetilde f(\bbu,\bar \bbu)$ is a continuous function on $\bc_\si^n$ 
and 
the strict transform $\tilde V$ of $V$ to $X$ in the toric chart $\bc_\si^n$ is given as follows: 
$$
\tilde V = \{ \widetilde f(\bbu,\bar \bbu) = 0 \} .
$$
We can also verify that 
$$
\{ \pi_\si^*f = 0 \} = \{ u_1 = 0 \} \cup \cdots \cup \{ u_k = 0 \} \cup \tilde V
$$
in the toric chart $\bc_\si^n$. %%%%

By \eqref{widetildeR0dots0}, we have 
$$
\begin{array}{ccl}
\tilde V(\tau)^* \cap \bc_\si^n & = & \tilde V \cap \{ u_1 = \cdots = u_k = 0,\ \text{and} \ u_{k+1} \neq 0, \cdots , u_n \neq 0 \} \\
                                          & = & \{ (\zero, \bbu')\ |\ \widetilde f_\De(\bbu',\bar\bbu') = 0,  \ u_{k+1} \neq 0, \cdots , u_n \neq 0 \} . 
\end{array}
$$
Note that $\tilde V(\tau)^* \cap \bc_\si^n$ is defined by $\widetilde f_\De(\bbu',\bar\bbu')$ only. 

Recall that 
\begin{equation}\label{f_De}
\pi_\si^*f_\De(\bbu,\bar \bbu) 
= \displaystyle \prod_{j=1}^k u_j^{\frac{r_j+p_j}2}\bar u_ j^{\frac{r_j-p_j}2} \, \times \widetilde f_\De(\bbu',\bar\bbu') .
\end{equation}
If $(\zero, \bx') \in \tilde V(\tau)^* \cap \bc_\si^n$, then we have 
$\widetilde f_\De(\bx',\bar\bx') = 0$, 
$$
\hat{\bx} := (\underbrace{1, \dots , 1}_{k} \, , \bx') \in \bc_\si^{*n} ,
$$
and $\pi_\si^*f_\De(\hat{\bx}) = \widetilde f_\De(\bx',\bar\bx') = 0$ by \eqref{f_De}. 
By the Newton non-degeneracy assumption on the compact face $\De$, 
$0$ is not a mixed critical \underline{value} of the face function $f_\De : \bc^{*n} \to \bc$. %%
Hence, $\hat{\bx}$ is not a mixed critical point of $\pi_\si^* f_\De$. 

For every $j \ (1 \leq j \leq k)$, we have 
$$
\begin{array}{cl}
\displaystyle 
\frac{\partial \pi_\si^*f_\De}{\partial u_j}  & = 
\displaystyle \frac{\partial}{\partial u_j} 
\left( \prod_{j=1}^k u_j^{\frac{r_j+p_j}2}\bar u_ j^{\frac{r_j-p_j}2} \right) \times \widetilde f_\De(\bbu',\bar\bbu') ,\\
\displaystyle 
\frac{\partial \pi_\si^*f_\De}{\partial \bar{u}_j}  & = 
\displaystyle \frac{\partial}{\partial \bar{u}_j} 
\left( \prod_{j=1}^k u_j^{\frac{r_j+p_j}2}\bar u_ j^{\frac{r_j-p_j}2} \right) \times \widetilde f_\De(\bbu',\bar\bbu') .
\end{array}
$$
Hence, for $\hat{\bx} = (\underbrace{1, \dots , 1}_{k}\, , \bx')$, we have 
$$
\begin{array}{cl}
\displaystyle 
\frac{\partial \pi_\si^*f_\De}{\partial u_j}(1, \dots , 1, \bx')  & = 
(\cdots \cdots) \times \widetilde f_\De(\bx',\bar\bx') = 0 ,\\
\displaystyle 
\frac{\partial \pi_\si^*f_\De}{\partial \bar{u}_j}(1, \dots , 1, \bx') & = 
(\cdots \cdots) \times \widetilde f_\De(\bx',\bar\bx') = 0 .
\end{array}
$$
Since $\hat{\bx} = (\underbrace{1, \dots , 1}_{k}\, , \bx')$ is a mixed regular point of $\pi_\si^* f_\De$, 
we see that 
$\bx'$ is a mixed regular point of the mixed polynomial $\widetilde f_\De(\bbu',\bar\bbu')$. 
By this result and the implicit function theorem, 
$\tilde V(\tau)^* \cap \bc_\si^n$ is a real analytic manifold of real dimension $2(n-k-1)$. 
This completes the proof of the assertion {\bf (i)}. 

\bigskip

We next prove the assertion {\bf (ii)}. 
Take an arbitrary point 
$$
\hat{\bbu}^0 \in \tilde V \cap \hat{\pi}^{-1}(\zero) = \tilde V \cap \bigcup_{P\gg 0}\hat{E}(P).
$$ 
By the {\sc Assumption}{\rm (*)}, 
there exists an $n$-dimensional cone $\sigma = (P_1, \dots , P_n)$ 
with strictly positive vertices $P_1, \dots , P_\kappa \ (1 \leq \kappa \leq n)$ 
for which the toric chart $\bc_\sigma^n$ contains $\hat{\bbu}^0$ 
and 
$$\hat{\bbu}^0 = \bbu_\si^0 = 
(\underbrace{0, \dots, 0}_{\kappa} \, ,\ u^0_{\kappa +1}, \dots, u^0_n),\ \ u^0_{\kappa +1}\neq 0, \dots ,u^0_n\neq 0
\ \ \ \text{in}\ \bc_\sigma^n .$$
We now use the symbol $k$ instead of $\kappa$, 
and consider the $k$-dimensional cone 
$$
\tau := (P_1, \dots , P_k) .
$$
Then we have 
$
\bbu_\si^0 \in \tilde V(\tau)^* \cap \bc_\si^n .
$ 
We can use the arguments and some results in the proof of {\bf (i)}. 
Then, we see that 
$(u^0_{k+1}, \dots, u^0_n)$ is a mixed regular point of  the function $\widetilde f_\De(\bbu',\bar\bbu')$. 
Put 
$$
N:= \sharp \{ (\nu, \mu) \ |\ c_{\nu, \mu}\neq 0,\ \nu + \mu \not\in \De \}, 
$$
and let the set 
$$
\{ (\nu_\ell, \mu_\ell) \ |\ \ell = 1, \dots , N \}
$$
coincide with the set
$
\{ (\nu, \mu) \ |\ \nu + \mu \not\in  \De \} .
$
Then we obtain another expression: 
$$
\widetilde R(\bbu,\bar\bbu) = 
\displaystyle \sum_{\ell = 1}^N 
c_\ell \left( \prod_{j=1}^k u_j^{a_j^\ell} \bar{u}_j^{b_j^\ell} \right) \left( \prod_{j=k+1}^n u_j^{a_j^\ell} \bar{u}_j^{b_j^\ell} \right) .
$$
Define a new polynomial with complex $N+(n-k)$ variables as follows: 
$$
\widetilde{\widetilde R}(v_1, \dots , v_N, \bbu') := 
\displaystyle \sum_{\ell = 1}^N c_\ell \, v_\ell \left( \prod_{j=k+1}^n u_j^{a_j^\ell} \bar{u}_j^{b_j^\ell} \right) .
$$
We write $\bv := (v_1, \dots , v_N)$, and define 
$$
\widetilde{\widetilde f}(\bv, \bbu') := \widetilde f_\De(\bbu',\bar\bbu') + \widetilde{\widetilde R}(\bv, \bbu') .
$$

\medskip

Now let $x_j$ and $y_j$, respectively, be the real part and the imaginary part of the variable $u_j$ 
for each $j \ (= k+1, \dots , n)$. 
Since $\widetilde{\widetilde R}(\zero, \bbu') = 0$ where $\zero = (\underbrace{0, \dots ,0}_N)$, 
the rank of the matrix 
$$
\begin{pmatrix}
\frac{\partial (\Re \widetilde{\widetilde f})}{\partial x_{k+1}} & \frac{\partial (\Re \widetilde{\widetilde f})}{\partial y_{k+1}} & \cdots \cdots & \frac{\partial (\Re \widetilde{\widetilde f})}{\partial x_n} & \frac{\partial (\Re \widetilde{\widetilde f})}{\partial y_n} \\
\frac{\partial (\Im \widetilde{\widetilde f})}{\partial x_{k+1}} & \frac{\partial (\Im \widetilde{\widetilde f})}{\partial y_{k+1}} & \cdots \cdots & \frac{\partial (\Im \widetilde{\widetilde f})}{\partial x_n} & \frac{\partial (\Im \widetilde{\widetilde f})}{\partial y_n}
\end{pmatrix}
$$
is $2$ at $(\zero, u^0_{k+1}, \dots, u^0_n)$. 
Especially, we see that $(\zero, u^0_{k+1}, \dots, u^0_n)$ is a mixed regular point of  the function $\widetilde{\widetilde f}(\bv, \bbu')$. 
Hence, by the implicit function theorem, 
there exists a real analytic map %%
$$F(\bv, x_{k+1}, y_{k+1}, \dots , \check{\xi} , \dots , \check{\eta} , \dots , x_n, y_n)$$
to $\br^2$ 
such that 
$$
\widetilde{\widetilde f}(\bv, \bbu') = 0 
\Longleftrightarrow 
\begin{pmatrix}
\xi \\
\eta
\end{pmatrix}
= 
F(\bv, x_{k+1}, y_{k+1}, \dots , \check{\xi} , \dots , \check{\eta} , \dots , x_n, y_n)
$$
in some neighborhood of the point $(\zero, u^0_{k+1}, \dots, u^0_n)$, where $\zero = (\underbrace{0, \dots ,0}_N)$. 
Here the notation $\check{\xi} , \dots , \check{\eta}$ means we remove two real variables $\xi$ and $\eta$. 
Hence, we have %%
$$
\begin{array}{cl}
 & \widetilde f(\bbu,\bar \bbu) = \widetilde{\widetilde f}(\left( \prod_{j=1}^k u_j^{a_j^1} \bar{u}_j^{b_j^1} \right), \dots , \left( \prod_{j=1}^k u_j^{a_j^N} \bar{u}_j^{b_j^N} \right), \bbu') = 0 \\
\Longleftrightarrow & 
\begin{pmatrix}
\xi \\
\eta
\end{pmatrix}
= 
F(\left( \prod_{j=1}^k u_j^{a_j^1} \bar{u}_j^{b_j^1} \right), \dots , \left( \prod_{j=1}^k u_j^{a_j^N} \bar{u}_j^{b_j^N} \right), x_{k+1}, y_{k+1}, \dots , \check{\xi} , \dots , \check{\eta} , \dots , x_n, y_n)
\end{array}
$$
in some neighborhood of the point $(\zero, u^0_{k+1}, \dots, u^0_n)$, 
where $\zero$ stands for $(\underbrace{0, \dots ,0}_k)$. 

\medskip

Suppose that $L(\Sigma^*) = \emptyset$. 
If $\tau$ is a cone of $\Sigma^*$ with only strictly positive vertices and 
$\tilde V(\tau)^* \cap \bc_\sigma^n \neq \emptyset$ for some $n$-dimensional cone $\sigma$ of $\Sigma^*$ with $\tau \prec \sigma$, 
then we see that 
$\{ P_j(\nu + \mu) - r_j\ |\ c_{\nu, \mu}\neq 0,\ \nu + \mu \not\in \De(P_j),\ j = 1, \dots, k \} = \emptyset$, where $k:= \dim \tau$. 
Hence, for every $(\nu, \mu)$ with $c_{\nu, \mu}\neq 0$, 
we have 
$\nu + \mu \in \De(P_j)$ for all $j = 1, \dots, k$, namely, 
$\nu + \mu \in \De$. 
Then we have $N=0$ and $\widetilde R(\bbu,\bar\bbu) \equiv 0$ by the argument in the proof of {\bf (i)}. 
Hence, $\tilde V(\tau)^* \cap \bc_\sigma^n$ is a real analytic manifold. 
Thus, $\tilde V$ is a real analytic manifold. %%

Next suppose that $L(\Sigma^*) \neq \emptyset$. 
Then there exists a $k$-dimensional cone $\tau$ of $\Sigma^*$ with only strictly positive vertices, 
$\{ P_j(\nu + \mu) - r_j\ |\ c_{\nu, \mu}\neq 0,\ \nu + \mu \not\in \De(P_j),\ j = 1, \dots, k \} \neq \emptyset$, and 
$\tilde V(\tau)^* \cap \bc_\sigma^n \neq \emptyset$ for some $n$-dimensional cone $\sigma$ of $\Sigma^*$ with $\tau \prec \sigma$. 
Then we have $N \geq 1$ and 
$\tilde V(\tau)^* \cap \bc_\sigma^n$ is a $C^{\Lambda(\tau) - 1}$-manifold by the argument in the proof of {\bf (i)}. 
Hence, $\tilde V$ is a $C^{\La - 1}$-manifold. 
Moreover, it is easy to see that 
$\tilde V$ is a real analytic manifold outside of $\tilde V \cap \hat{\pi}^{-1}(\zero)$. %%

This completes the proof of the assertion {\bf (ii)}. 
\end{proof}

\medskip

\begin{corollary}[mixed homogeneous polynomial case,\ $n = 2$]
Let us consider the case $n = 2$. %%
Let $(f,\zero)$ be the germ at $\zero \in \bc^2$ of 
a convenient and Newton non-degenerate strongly polar non-negative mixed weighted homogeneous polynomial 
$f(\bbz, \bar{\bbz}) = \sum_{\nu, \mu}c_{\nu, \mu} \bbz^\nu \bar{\bbz}^\mu $
of radial degree $r$ with respect to the weight vector $P := {}^t(1,1)$, and hence, 
$f(\bbz, \bar{\bbz})$ is a mixed \underline{homogeneous} polynomial. 
We assume $f(\zero)=0$. 
Then, all the compact faces of the Newton polyhedron $\Gamma_+(f)$ 
are the $1$-dimensional face $\De(P)$ and the $0$-dimensional faces $\De(P+E_1),\ \De(P+E_2)$. 
By the assumption, we have $\pdeg_Pf_P = \pdeg_Pf \geq 0$. 
Then it is easy to check (Recall Lemma \ref{E-j}) that 
$\rdeg_{P+E_i}f_{P+E_i} = \rdeg_Pf_{P+E_i} = \rdeg_Pf_P = r$ and 
$\pdeg_{P+E_i}f_{P+E_i} = \pdeg_Pf_{P+E_i} = \pdeg_Pf_P \geq 0$ for $i=1,2$. 
Hence, we can say that 
$(f,\zero)$ is of strongly polar non-negative mixed weighted homogeneous face type (Definition \ref{strongly-mixed-wh-face-type}). 
In this case the (radial) dual Newton diagram $\Gamma^*(f)$ 
gives already a regular simplicial subdivision of $N^+_\br$ which is admissible for $f$. 
Moreover, $\Gamma^*(f)$ is obviously convenient in the sense of Definition \ref{convenient-subdivision-f}. 
Consider the toric modification $\hat{\pi} : X \to \bc^n$ associated with $\Gamma^*(f)$. 
This is nothing but the ``usual blowing-up". 
Let $(V,\zero)$ be the germ of the mixed hypersurface $V := f^{-1}(0)$ at $\zero \in \bc^2$. 
We assume the {\sc Assumption}{\rm (*)}. 
Since $\{ P(\nu + \mu) - r\ |\ c_{\nu, \mu}\neq 0,\ \nu + \mu \not\in \De(P) \} = \emptyset$, 
we have $L(\Sigma^*) = L(\Gamma^*(f)) = \emptyset$. 
By Theorem \ref{theorem11-improved} {\bf (ii)}, we conclude that 
the strict transform germ $\tilde V$ at $\tilde V \cap \hat{\pi}^{-1}(\zero)$ 
is a real analytic manifold. 
\end{corollary}

\medskip

\section{The mixed polynomial $z_1^4 \bar{z}_1^2 + z_1^a \bar{z}_1^{4-a} z_2^b \bar{z}_2^{3-b} + \bar{z}_1^2 z_2^3 + z_2^3 \bar{z}_2^3$} \label{section-exam4-3-mixed}

As a typical example of Theorem \ref{theorem11-improved}, 
let us consider the following mixed polynomial of 
strongly polar non-negative mixed weighted homogeneous face type (see Definition \ref{strongly-mixed-wh-face-type})
$$f(\bbz, \bar{\bbz}) := (z_1^4 + z_2^3)\overline{(z_1^2 + z_2^3)}.$$

Example 4.3 %%%%
 of \cite{Oka2018}, p.73 (see also \cite{Saito-Takashimizu2021-2}) 
describes a good resolution of the germ at $\zero$ of the holomorphic polynomial 
$$(z_1^4 + z_2^3)(z_1^2 + z_2^3).$$
Since $(z_1^4 + z_2^3)\overline{(z_1^2 + z_2^3)}$ and $(z_1^4 + z_2^3)(z_1^2 + z_2^3)$ have 
the same zero set and the same Newton polyhedron, 
we must be able to resolve the isolated singular point $\zero$ of 
the germ $((z_1^4 + z_2^3)\overline{(z_1^2 + z_2^3)}, \zero)$ 
by the same toric modification. 
We have 
\begin{equation}\label{exam4-3-mixed}
f(\bbz, \bar{\bbz}) = (z_1^4 + z_2^3)\overline{(z_1^2 + z_2^3)} = z_1^4\bar{z}_1^2 + z_1^4\bar{z}_2^3 + \bar{z}_1^2z_2^3 + z_2^3\bar{z}_2^3 .
\end{equation}

\medskip

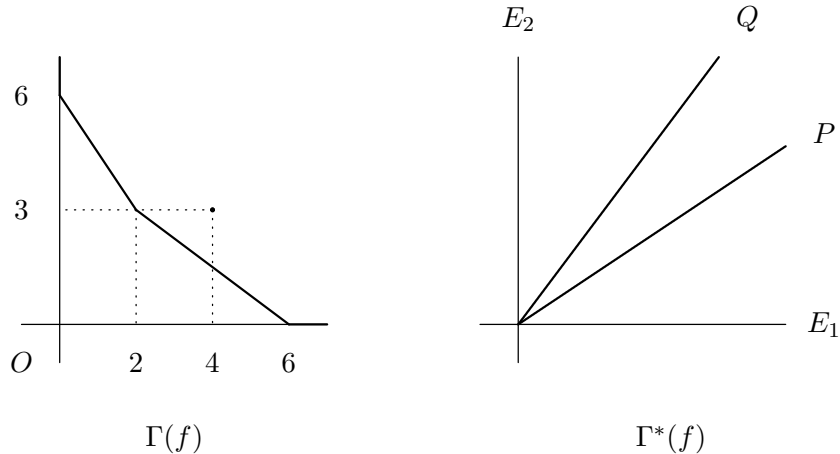
\begin{figure}[H]
\begin{center}
%% \input{exm.4.3_newton_diagram}
%WinTpicVersion4.32a
{\unitlength 0.1in%
\begin{picture}(42.1000,22.0000)(3.0000,-25.3500)%
% LINE 2 0 3 0 Black White  
% 8 600 600 600 2200 400 2000 2000 2000 3000 600 3000 2200 2800 2000 4400 2000
% 
\special{pn 8}%
\special{pa 600 600}%
\special{pa 600 2200}%
\special{fp}%
\special{pa 400 2000}%
\special{pa 2000 2000}%
\special{fp}%
\special{pa 3000 600}%
\special{pa 3000 2200}%
\special{fp}%
\special{pa 2800 2000}%
\special{pa 4400 2000}%
\special{fp}%
% LINE 1 0 3 0 Black White  
% 8 600 800 1000 1400 1000 1400 1800 2000 1800 2000 2000 2000 600 800 600 600
% 
\special{pn 13}%
\special{pa 600 800}%
\special{pa 1000 1400}%
\special{fp}%
\special{pa 1000 1400}%
\special{pa 1800 2000}%
\special{fp}%
\special{pa 1800 2000}%
\special{pa 2000 2000}%
\special{fp}%
\special{pa 600 800}%
\special{pa 600 600}%
\special{fp}%
% FUNC 1 0 3 0 Black White  
% 9 3000 600 4400 2000 3000 2000 3200 2000 3000 1800 3000 600 4400 2000 0 2 0 0
% 2x/3
\special{pn 13}%
\special{pa 3000 2000}%
\special{pa 3005 1997}%
\special{pa 3010 1993}%
\special{pa 3020 1987}%
\special{pa 3025 1983}%
\special{pa 3035 1977}%
\special{pa 3040 1973}%
\special{pa 3050 1967}%
\special{pa 3055 1963}%
\special{pa 3065 1957}%
\special{pa 3070 1953}%
\special{pa 3080 1947}%
\special{pa 3085 1943}%
\special{pa 3095 1937}%
\special{pa 3100 1933}%
\special{pa 3110 1927}%
\special{pa 3115 1923}%
\special{pa 3125 1917}%
\special{pa 3130 1913}%
\special{pa 3140 1907}%
\special{pa 3145 1903}%
\special{pa 3155 1897}%
\special{pa 3160 1893}%
\special{pa 3170 1887}%
\special{pa 3175 1883}%
\special{pa 3185 1877}%
\special{pa 3190 1873}%
\special{pa 3200 1867}%
\special{pa 3205 1863}%
\special{pa 3215 1857}%
\special{pa 3220 1853}%
\special{pa 3230 1847}%
\special{pa 3235 1843}%
\special{pa 3245 1837}%
\special{pa 3250 1833}%
\special{pa 3260 1827}%
\special{pa 3265 1823}%
\special{pa 3275 1817}%
\special{pa 3280 1813}%
\special{pa 3290 1807}%
\special{pa 3295 1803}%
\special{pa 3305 1797}%
\special{pa 3310 1793}%
\special{pa 3320 1787}%
\special{pa 3325 1783}%
\special{pa 3335 1777}%
\special{pa 3340 1773}%
\special{pa 3350 1767}%
\special{pa 3355 1763}%
\special{pa 3365 1757}%
\special{pa 3370 1753}%
\special{pa 3380 1747}%
\special{pa 3385 1743}%
\special{pa 3395 1737}%
\special{pa 3400 1733}%
\special{pa 3410 1727}%
\special{pa 3415 1723}%
\special{pa 3425 1717}%
\special{pa 3430 1713}%
\special{pa 3440 1707}%
\special{pa 3445 1703}%
\special{pa 3455 1697}%
\special{pa 3460 1693}%
\special{pa 3470 1687}%
\special{pa 3475 1683}%
\special{pa 3485 1677}%
\special{pa 3490 1673}%
\special{pa 3500 1667}%
\special{pa 3505 1663}%
\special{pa 3515 1657}%
\special{pa 3520 1653}%
\special{pa 3530 1647}%
\special{pa 3535 1643}%
\special{pa 3545 1637}%
\special{pa 3550 1633}%
\special{pa 3560 1627}%
\special{pa 3565 1623}%
\special{pa 3575 1617}%
\special{pa 3580 1613}%
\special{pa 3590 1607}%
\special{pa 3595 1603}%
\special{pa 3605 1597}%
\special{pa 3610 1593}%
\special{pa 3620 1587}%
\special{pa 3625 1583}%
\special{pa 3635 1577}%
\special{pa 3640 1573}%
\special{pa 3650 1567}%
\special{pa 3655 1563}%
\special{pa 3665 1557}%
\special{pa 3670 1553}%
\special{pa 3680 1547}%
\special{pa 3685 1543}%
\special{pa 3695 1537}%
\special{pa 3700 1533}%
\special{pa 3710 1527}%
\special{pa 3715 1523}%
\special{pa 3725 1517}%
\special{pa 3730 1513}%
\special{pa 3740 1507}%
\special{pa 3745 1503}%
\special{pa 3755 1497}%
\special{pa 3760 1493}%
\special{pa 3770 1487}%
\special{pa 3775 1483}%
\special{pa 3785 1477}%
\special{pa 3790 1473}%
\special{pa 3800 1467}%
\special{pa 3805 1463}%
\special{pa 3815 1457}%
\special{pa 3820 1453}%
\special{pa 3830 1447}%
\special{pa 3835 1443}%
\special{pa 3845 1437}%
\special{pa 3850 1433}%
\special{pa 3860 1427}%
\special{pa 3865 1423}%
\special{pa 3875 1417}%
\special{pa 3880 1413}%
\special{pa 3890 1407}%
\special{pa 3895 1403}%
\special{pa 3905 1397}%
\special{pa 3910 1393}%
\special{pa 3920 1387}%
\special{pa 3925 1383}%
\special{pa 3935 1377}%
\special{pa 3940 1373}%
\special{pa 3950 1367}%
\special{pa 3955 1363}%
\special{pa 3965 1357}%
\special{pa 3970 1353}%
\special{pa 3980 1347}%
\special{pa 3985 1343}%
\special{pa 3995 1337}%
\special{pa 4000 1333}%
\special{pa 4010 1327}%
\special{pa 4015 1323}%
\special{pa 4025 1317}%
\special{pa 4030 1313}%
\special{pa 4040 1307}%
\special{pa 4045 1303}%
\special{pa 4055 1297}%
\special{pa 4060 1293}%
\special{pa 4070 1287}%
\special{pa 4075 1283}%
\special{pa 4085 1277}%
\special{pa 4090 1273}%
\special{pa 4100 1267}%
\special{pa 4105 1263}%
\special{pa 4115 1257}%
\special{pa 4120 1253}%
\special{pa 4130 1247}%
\special{pa 4135 1243}%
\special{pa 4145 1237}%
\special{pa 4150 1233}%
\special{pa 4160 1227}%
\special{pa 4165 1223}%
\special{pa 4175 1217}%
\special{pa 4180 1213}%
\special{pa 4190 1207}%
\special{pa 4195 1203}%
\special{pa 4205 1197}%
\special{pa 4210 1193}%
\special{pa 4220 1187}%
\special{pa 4225 1183}%
\special{pa 4235 1177}%
\special{pa 4240 1173}%
\special{pa 4250 1167}%
\special{pa 4255 1163}%
\special{pa 4265 1157}%
\special{pa 4270 1153}%
\special{pa 4280 1147}%
\special{pa 4285 1143}%
\special{pa 4295 1137}%
\special{pa 4300 1133}%
\special{pa 4310 1127}%
\special{pa 4315 1123}%
\special{pa 4325 1117}%
\special{pa 4330 1113}%
\special{pa 4340 1107}%
\special{pa 4345 1103}%
\special{pa 4355 1097}%
\special{pa 4360 1093}%
\special{pa 4370 1087}%
\special{pa 4375 1083}%
\special{pa 4385 1077}%
\special{pa 4390 1073}%
\special{pa 4400 1067}%
\special{fp}%
% FUNC 1 0 3 0 Black White  
% 10 3000 600 4400 2000 3000 2000 3200 2000 3000 1800 3000 600 4400 2000 0 2 0 0 0 0
% 4x/3
\special{pn 13}%
\special{pa 3000 2000}%
\special{pa 3005 1993}%
\special{pa 3010 1987}%
\special{pa 3020 1973}%
\special{pa 3025 1967}%
\special{pa 3035 1953}%
\special{pa 3040 1947}%
\special{pa 3050 1933}%
\special{pa 3055 1927}%
\special{pa 3065 1913}%
\special{pa 3070 1907}%
\special{pa 3080 1893}%
\special{pa 3085 1887}%
\special{pa 3095 1873}%
\special{pa 3100 1867}%
\special{pa 3110 1853}%
\special{pa 3115 1847}%
\special{pa 3125 1833}%
\special{pa 3130 1827}%
\special{pa 3140 1813}%
\special{pa 3145 1807}%
\special{pa 3155 1793}%
\special{pa 3160 1787}%
\special{pa 3170 1773}%
\special{pa 3175 1767}%
\special{pa 3185 1753}%
\special{pa 3190 1747}%
\special{pa 3200 1733}%
\special{pa 3205 1727}%
\special{pa 3215 1713}%
\special{pa 3220 1707}%
\special{pa 3230 1693}%
\special{pa 3235 1687}%
\special{pa 3245 1673}%
\special{pa 3250 1667}%
\special{pa 3260 1653}%
\special{pa 3265 1647}%
\special{pa 3275 1633}%
\special{pa 3280 1627}%
\special{pa 3290 1613}%
\special{pa 3295 1607}%
\special{pa 3305 1593}%
\special{pa 3310 1587}%
\special{pa 3320 1573}%
\special{pa 3325 1567}%
\special{pa 3335 1553}%
\special{pa 3340 1547}%
\special{pa 3350 1533}%
\special{pa 3355 1527}%
\special{pa 3365 1513}%
\special{pa 3370 1507}%
\special{pa 3380 1493}%
\special{pa 3385 1487}%
\special{pa 3395 1473}%
\special{pa 3400 1467}%
\special{pa 3410 1453}%
\special{pa 3415 1447}%
\special{pa 3425 1433}%
\special{pa 3430 1427}%
\special{pa 3440 1413}%
\special{pa 3445 1407}%
\special{pa 3455 1393}%
\special{pa 3460 1387}%
\special{pa 3470 1373}%
\special{pa 3475 1367}%
\special{pa 3485 1353}%
\special{pa 3490 1347}%
\special{pa 3500 1333}%
\special{pa 3505 1327}%
\special{pa 3515 1313}%
\special{pa 3520 1307}%
\special{pa 3530 1293}%
\special{pa 3535 1287}%
\special{pa 3545 1273}%
\special{pa 3550 1267}%
\special{pa 3560 1253}%
\special{pa 3565 1247}%
\special{pa 3575 1233}%
\special{pa 3580 1227}%
\special{pa 3590 1213}%
\special{pa 3595 1207}%
\special{pa 3605 1193}%
\special{pa 3610 1187}%
\special{pa 3620 1173}%
\special{pa 3625 1167}%
\special{pa 3635 1153}%
\special{pa 3640 1147}%
\special{pa 3650 1133}%
\special{pa 3655 1127}%
\special{pa 3665 1113}%
\special{pa 3670 1107}%
\special{pa 3680 1093}%
\special{pa 3685 1087}%
\special{pa 3695 1073}%
\special{pa 3700 1067}%
\special{pa 3710 1053}%
\special{pa 3715 1047}%
\special{pa 3725 1033}%
\special{pa 3730 1027}%
\special{pa 3740 1013}%
\special{pa 3745 1007}%
\special{pa 3755 993}%
\special{pa 3760 987}%
\special{pa 3770 973}%
\special{pa 3775 967}%
\special{pa 3785 953}%
\special{pa 3790 947}%
\special{pa 3800 933}%
\special{pa 3805 927}%
\special{pa 3815 913}%
\special{pa 3820 907}%
\special{pa 3830 893}%
\special{pa 3835 887}%
\special{pa 3845 873}%
\special{pa 3850 867}%
\special{pa 3860 853}%
\special{pa 3865 847}%
\special{pa 3875 833}%
\special{pa 3880 827}%
\special{pa 3890 813}%
\special{pa 3895 807}%
\special{pa 3905 793}%
\special{pa 3910 787}%
\special{pa 3920 773}%
\special{pa 3925 767}%
\special{pa 3935 753}%
\special{pa 3940 747}%
\special{pa 3950 733}%
\special{pa 3955 727}%
\special{pa 3965 713}%
\special{pa 3970 707}%
\special{pa 3980 693}%
\special{pa 3985 687}%
\special{pa 3995 673}%
\special{pa 4000 667}%
\special{pa 4010 653}%
\special{pa 4015 647}%
\special{pa 4025 633}%
\special{pa 4030 627}%
\special{pa 4040 613}%
\special{pa 4045 607}%
\special{pa 4050 600}%
\special{fp}%
% LINE 2 2 3 0 Black White  
% 4 1000 1400 600 1400 1000 1400 1000 2000
% 
\special{pn 8}%
\special{pa 1000 1400}%
\special{pa 600 1400}%
\special{dt 0.045}%
\special{pa 1000 1400}%
\special{pa 1000 2000}%
\special{dt 0.045}%
% STR 2 0 3 0 Black White  
% 4 400 2100 400 2200 5 0 0 0
% $O$
\put(4.0000,-22.0000){\makebox(0,0){$O$}}%
% STR 2 0 3 0 Black White  
% 4 400 700 400 800 5 0 0 0
% $6$
\put(4.0000,-8.0000){\makebox(0,0){$6$}}%
% STR 2 0 3 0 Black White  
% 4 400 1300 400 1400 5 0 0 0
% $3$
\put(4.0000,-14.0000){\makebox(0,0){$3$}}%
% STR 2 0 3 0 Black White  
% 4 1000 2100 1000 2200 5 0 0 0
% $2$
\put(10.0000,-22.0000){\makebox(0,0){$2$}}%
% DOT 1 2 3 0 Black White  
% 2 1400 1400 1400 1400
% 
\special{pn 4}%
\special{sh 1}%
\special{ar 1400 1400 10 10 0 6.2831853}%
\special{sh 1}%
\special{ar 1400 1400 10 10 0 6.2831853}%
% LINE 2 2 3 0 Black White  
% 4 1400 1400 1000 1400 1400 1400 1400 2000
% 
\special{pn 8}%
\special{pa 1400 1400}%
\special{pa 1000 1400}%
\special{dt 0.045}%
\special{pa 1400 1400}%
\special{pa 1400 2000}%
\special{dt 0.045}%
% STR 2 0 3 0 Black White  
% 4 1400 2100 1400 2200 5 0 0 0
% $4$
\put(14.0000,-22.0000){\makebox(0,0){$4$}}%
% STR 2 0 3 0 Black White  
% 4 1800 2100 1800 2200 5 0 0 0
% $6$
\put(18.0000,-22.0000){\makebox(0,0){$6$}}%
% STR 2 0 3 0 Black White  
% 4 4600 1900 4600 2000 5 0 0 0
% $E_{1}$
\put(46.0000,-20.0000){\makebox(0,0){$E_{1}$}}%
% STR 2 0 3 0 Black White  
% 4 4600 900 4600 1000 5 0 0 0
% $P$
\put(46.0000,-10.0000){\makebox(0,0){$P$}}%
% STR 2 0 3 0 Black White  
% 4 4200 300 4200 400 5 0 0 0
% $Q$
\put(42.0000,-4.0000){\makebox(0,0){$Q$}}%
% STR 2 0 3 0 Black White  
% 4 3000 300 3000 400 5 0 0 0
% $E_{2}$
\put(30.0000,-4.0000){\makebox(0,0){$E_{2}$}}%
% STR 2 0 3 0 Black White  
% 4 1200 2500 1200 2600 5 0 0 0
% $\Gamma(f)$
\put(12.0000,-26.0000){\makebox(0,0){$\Gamma(f)$}}%
% STR 2 0 3 0 Black White  
% 4 3800 2500 3800 2600 5 0 0 0
% $\Gamma^{*}(f)$
\put(38.0000,-26.0000){\makebox(0,0){$\Gamma^{*}(f)$}}%
\end{picture}}%
\end{center}

\vspace{7mm}

\caption{The radial Newton polyhedron $\Gamma_+(f)$ and the dual Newton diagram $\Gamma^*(f)$.}
\label{exm4_3_newton_diagram}
\end{figure}

The radial Newton polyhedron $\Gamma_+(f)$ of $(f,\zero)$ and its dual Newton diagram $\Gamma^*(f)$ 
are as in Figure \ref{exm4_3_newton_diagram}, 
where we set 
$$E_1 = {}^t(1,0),\ P= {}^t(3,2),\ Q= {}^t(3,4),\ E_2 = {}^t(0,1).$$

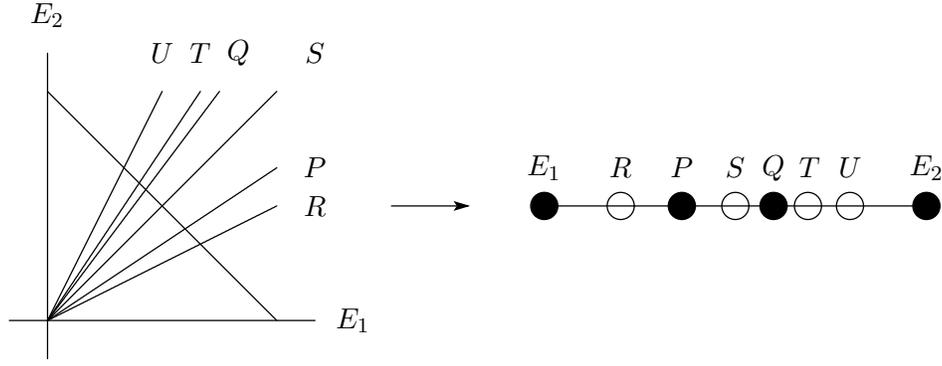
\begin{figure}[H]
\begin{center}
%%\input{exm.4.3_regular_subdivision}
%WinTpicVersion4.32a
{\unitlength 0.1in%
\begin{picture}(48.7100,18.6500)(12.0000,-26.0000)%
% LINE 2 0 3 0 Black White  
% 6 1400 1000 1400 2600 1200 2400 2800 2400 4000 1800 6000 1800
% 
\special{pn 8}%
\special{pa 1400 1000}%
\special{pa 1400 2600}%
\special{fp}%
\special{pa 1200 2400}%
\special{pa 2800 2400}%
\special{fp}%
\special{pa 4000 1800}%
\special{pa 6000 1800}%
\special{fp}%
% VECTOR 2 0 3 0 Black White  
% 2 3200 1800 3600 1800
% 
\special{pn 8}%
\special{pa 3200 1800}%
\special{pa 3600 1800}%
\special{fp}%
\special{sh 1}%
\special{pa 3600 1800}%
\special{pa 3533 1780}%
\special{pa 3547 1800}%
\special{pa 3533 1820}%
\special{pa 3600 1800}%
\special{fp}%
% FUNC 2 0 3 0 Black White  
% 9 1400 1200 2600 2400 1400 2400 1600 2400 1400 2200 1400 1200 2600 2400 0 2 0 0
% x/2
\special{pn 8}%
\special{pa 1400 2400}%
\special{pa 1405 2398}%
\special{pa 1415 2392}%
\special{pa 1425 2388}%
\special{pa 1435 2382}%
\special{pa 1445 2378}%
\special{pa 1455 2372}%
\special{pa 1465 2368}%
\special{pa 1475 2362}%
\special{pa 1485 2358}%
\special{pa 1495 2352}%
\special{pa 1505 2348}%
\special{pa 1515 2342}%
\special{pa 1525 2338}%
\special{pa 1535 2332}%
\special{pa 1545 2328}%
\special{pa 1555 2322}%
\special{pa 1565 2318}%
\special{pa 1575 2312}%
\special{pa 1585 2308}%
\special{pa 1595 2302}%
\special{pa 1605 2298}%
\special{pa 1615 2292}%
\special{pa 1625 2288}%
\special{pa 1635 2282}%
\special{pa 1645 2278}%
\special{pa 1655 2272}%
\special{pa 1665 2268}%
\special{pa 1675 2262}%
\special{pa 1685 2258}%
\special{pa 1695 2252}%
\special{pa 1705 2248}%
\special{pa 1715 2242}%
\special{pa 1725 2238}%
\special{pa 1735 2232}%
\special{pa 1745 2228}%
\special{pa 1755 2222}%
\special{pa 1765 2218}%
\special{pa 1775 2212}%
\special{pa 1785 2208}%
\special{pa 1795 2202}%
\special{pa 1805 2198}%
\special{pa 1815 2192}%
\special{pa 1825 2188}%
\special{pa 1835 2182}%
\special{pa 1845 2178}%
\special{pa 1855 2172}%
\special{pa 1865 2168}%
\special{pa 1875 2162}%
\special{pa 1885 2158}%
\special{pa 1895 2152}%
\special{pa 1905 2148}%
\special{pa 1915 2142}%
\special{pa 1925 2138}%
\special{pa 1935 2132}%
\special{pa 1945 2128}%
\special{pa 1955 2122}%
\special{pa 1965 2118}%
\special{pa 1975 2112}%
\special{pa 1985 2108}%
\special{pa 1995 2102}%
\special{pa 2005 2098}%
\special{pa 2015 2092}%
\special{pa 2025 2088}%
\special{pa 2035 2082}%
\special{pa 2045 2078}%
\special{pa 2055 2072}%
\special{pa 2065 2068}%
\special{pa 2075 2062}%
\special{pa 2085 2058}%
\special{pa 2095 2052}%
\special{pa 2105 2048}%
\special{pa 2115 2042}%
\special{pa 2125 2038}%
\special{pa 2135 2032}%
\special{pa 2145 2028}%
\special{pa 2155 2022}%
\special{pa 2165 2018}%
\special{pa 2175 2012}%
\special{pa 2185 2008}%
\special{pa 2195 2002}%
\special{pa 2205 1998}%
\special{pa 2215 1992}%
\special{pa 2225 1988}%
\special{pa 2235 1982}%
\special{pa 2245 1978}%
\special{pa 2255 1972}%
\special{pa 2265 1968}%
\special{pa 2275 1962}%
\special{pa 2285 1958}%
\special{pa 2295 1952}%
\special{pa 2305 1948}%
\special{pa 2315 1942}%
\special{pa 2325 1938}%
\special{pa 2335 1932}%
\special{pa 2345 1928}%
\special{pa 2355 1922}%
\special{pa 2365 1918}%
\special{pa 2375 1912}%
\special{pa 2385 1908}%
\special{pa 2395 1902}%
\special{pa 2405 1898}%
\special{pa 2415 1892}%
\special{pa 2425 1888}%
\special{pa 2435 1882}%
\special{pa 2445 1878}%
\special{pa 2455 1872}%
\special{pa 2465 1868}%
\special{pa 2475 1862}%
\special{pa 2485 1858}%
\special{pa 2495 1852}%
\special{pa 2505 1848}%
\special{pa 2515 1842}%
\special{pa 2525 1838}%
\special{pa 2535 1832}%
\special{pa 2545 1828}%
\special{pa 2555 1822}%
\special{pa 2565 1818}%
\special{pa 2575 1812}%
\special{pa 2585 1808}%
\special{pa 2595 1802}%
\special{pa 2600 1800}%
\special{fp}%
% FUNC 2 0 3 0 Black White  
% 10 1400 1200 2600 2400 1400 2400 1600 2400 1400 2200 1400 1200 2600 2400 0 2 0 0 0 0
% 2x/3
\special{pn 8}%
\special{pa 1400 2400}%
\special{pa 1405 2397}%
\special{pa 1410 2393}%
\special{pa 1420 2387}%
\special{pa 1425 2383}%
\special{pa 1435 2377}%
\special{pa 1440 2373}%
\special{pa 1450 2367}%
\special{pa 1455 2363}%
\special{pa 1465 2357}%
\special{pa 1470 2353}%
\special{pa 1480 2347}%
\special{pa 1485 2343}%
\special{pa 1495 2337}%
\special{pa 1500 2333}%
\special{pa 1510 2327}%
\special{pa 1515 2323}%
\special{pa 1525 2317}%
\special{pa 1530 2313}%
\special{pa 1540 2307}%
\special{pa 1545 2303}%
\special{pa 1555 2297}%
\special{pa 1560 2293}%
\special{pa 1570 2287}%
\special{pa 1575 2283}%
\special{pa 1585 2277}%
\special{pa 1590 2273}%
\special{pa 1600 2267}%
\special{pa 1605 2263}%
\special{pa 1615 2257}%
\special{pa 1620 2253}%
\special{pa 1630 2247}%
\special{pa 1635 2243}%
\special{pa 1645 2237}%
\special{pa 1650 2233}%
\special{pa 1660 2227}%
\special{pa 1665 2223}%
\special{pa 1675 2217}%
\special{pa 1680 2213}%
\special{pa 1690 2207}%
\special{pa 1695 2203}%
\special{pa 1705 2197}%
\special{pa 1710 2193}%
\special{pa 1720 2187}%
\special{pa 1725 2183}%
\special{pa 1735 2177}%
\special{pa 1740 2173}%
\special{pa 1750 2167}%
\special{pa 1755 2163}%
\special{pa 1765 2157}%
\special{pa 1770 2153}%
\special{pa 1780 2147}%
\special{pa 1785 2143}%
\special{pa 1795 2137}%
\special{pa 1800 2133}%
\special{pa 1810 2127}%
\special{pa 1815 2123}%
\special{pa 1825 2117}%
\special{pa 1830 2113}%
\special{pa 1840 2107}%
\special{pa 1845 2103}%
\special{pa 1855 2097}%
\special{pa 1860 2093}%
\special{pa 1870 2087}%
\special{pa 1875 2083}%
\special{pa 1885 2077}%
\special{pa 1890 2073}%
\special{pa 1900 2067}%
\special{pa 1905 2063}%
\special{pa 1915 2057}%
\special{pa 1920 2053}%
\special{pa 1930 2047}%
\special{pa 1935 2043}%
\special{pa 1945 2037}%
\special{pa 1950 2033}%
\special{pa 1960 2027}%
\special{pa 1965 2023}%
\special{pa 1975 2017}%
\special{pa 1980 2013}%
\special{pa 1990 2007}%
\special{pa 1995 2003}%
\special{pa 2005 1997}%
\special{pa 2010 1993}%
\special{pa 2020 1987}%
\special{pa 2025 1983}%
\special{pa 2035 1977}%
\special{pa 2040 1973}%
\special{pa 2050 1967}%
\special{pa 2055 1963}%
\special{pa 2065 1957}%
\special{pa 2070 1953}%
\special{pa 2080 1947}%
\special{pa 2085 1943}%
\special{pa 2095 1937}%
\special{pa 2100 1933}%
\special{pa 2110 1927}%
\special{pa 2115 1923}%
\special{pa 2125 1917}%
\special{pa 2130 1913}%
\special{pa 2140 1907}%
\special{pa 2145 1903}%
\special{pa 2155 1897}%
\special{pa 2160 1893}%
\special{pa 2170 1887}%
\special{pa 2175 1883}%
\special{pa 2185 1877}%
\special{pa 2190 1873}%
\special{pa 2200 1867}%
\special{pa 2205 1863}%
\special{pa 2215 1857}%
\special{pa 2220 1853}%
\special{pa 2230 1847}%
\special{pa 2235 1843}%
\special{pa 2245 1837}%
\special{pa 2250 1833}%
\special{pa 2260 1827}%
\special{pa 2265 1823}%
\special{pa 2275 1817}%
\special{pa 2280 1813}%
\special{pa 2290 1807}%
\special{pa 2295 1803}%
\special{pa 2305 1797}%
\special{pa 2310 1793}%
\special{pa 2320 1787}%
\special{pa 2325 1783}%
\special{pa 2335 1777}%
\special{pa 2340 1773}%
\special{pa 2350 1767}%
\special{pa 2355 1763}%
\special{pa 2365 1757}%
\special{pa 2370 1753}%
\special{pa 2380 1747}%
\special{pa 2385 1743}%
\special{pa 2395 1737}%
\special{pa 2400 1733}%
\special{pa 2410 1727}%
\special{pa 2415 1723}%
\special{pa 2425 1717}%
\special{pa 2430 1713}%
\special{pa 2440 1707}%
\special{pa 2445 1703}%
\special{pa 2455 1697}%
\special{pa 2460 1693}%
\special{pa 2470 1687}%
\special{pa 2475 1683}%
\special{pa 2485 1677}%
\special{pa 2490 1673}%
\special{pa 2500 1667}%
\special{pa 2505 1663}%
\special{pa 2515 1657}%
\special{pa 2520 1653}%
\special{pa 2530 1647}%
\special{pa 2535 1643}%
\special{pa 2545 1637}%
\special{pa 2550 1633}%
\special{pa 2560 1627}%
\special{pa 2565 1623}%
\special{pa 2575 1617}%
\special{pa 2580 1613}%
\special{pa 2590 1607}%
\special{pa 2595 1603}%
\special{pa 2600 1600}%
\special{fp}%
% FUNC 2 0 3 0 Black White  
% 10 1400 1200 2600 2400 1400 2400 1600 2400 1400 2200 1400 1200 2600 2400 0 2 0 0 0 0
% x
\special{pn 8}%
\special{pa 1400 2400}%
\special{pa 2600 1200}%
\special{fp}%
% FUNC 2 0 3 0 Black White  
% 10 1400 1200 2600 2400 1400 2400 1600 2400 1400 2200 1400 1200 2600 2400 0 2 0 0 0 0
% 4x/3
\special{pn 8}%
\special{pa 1400 2400}%
\special{pa 1405 2393}%
\special{pa 1410 2387}%
\special{pa 1420 2373}%
\special{pa 1425 2367}%
\special{pa 1435 2353}%
\special{pa 1440 2347}%
\special{pa 1450 2333}%
\special{pa 1455 2327}%
\special{pa 1465 2313}%
\special{pa 1470 2307}%
\special{pa 1480 2293}%
\special{pa 1485 2287}%
\special{pa 1495 2273}%
\special{pa 1500 2267}%
\special{pa 1510 2253}%
\special{pa 1515 2247}%
\special{pa 1525 2233}%
\special{pa 1530 2227}%
\special{pa 1540 2213}%
\special{pa 1545 2207}%
\special{pa 1555 2193}%
\special{pa 1560 2187}%
\special{pa 1570 2173}%
\special{pa 1575 2167}%
\special{pa 1585 2153}%
\special{pa 1590 2147}%
\special{pa 1600 2133}%
\special{pa 1605 2127}%
\special{pa 1615 2113}%
\special{pa 1620 2107}%
\special{pa 1630 2093}%
\special{pa 1635 2087}%
\special{pa 1645 2073}%
\special{pa 1650 2067}%
\special{pa 1660 2053}%
\special{pa 1665 2047}%
\special{pa 1675 2033}%
\special{pa 1680 2027}%
\special{pa 1690 2013}%
\special{pa 1695 2007}%
\special{pa 1705 1993}%
\special{pa 1710 1987}%
\special{pa 1720 1973}%
\special{pa 1725 1967}%
\special{pa 1735 1953}%
\special{pa 1740 1947}%
\special{pa 1750 1933}%
\special{pa 1755 1927}%
\special{pa 1765 1913}%
\special{pa 1770 1907}%
\special{pa 1780 1893}%
\special{pa 1785 1887}%
\special{pa 1795 1873}%
\special{pa 1800 1867}%
\special{pa 1810 1853}%
\special{pa 1815 1847}%
\special{pa 1825 1833}%
\special{pa 1830 1827}%
\special{pa 1840 1813}%
\special{pa 1845 1807}%
\special{pa 1855 1793}%
\special{pa 1860 1787}%
\special{pa 1870 1773}%
\special{pa 1875 1767}%
\special{pa 1885 1753}%
\special{pa 1890 1747}%
\special{pa 1900 1733}%
\special{pa 1905 1727}%
\special{pa 1915 1713}%
\special{pa 1920 1707}%
\special{pa 1930 1693}%
\special{pa 1935 1687}%
\special{pa 1945 1673}%
\special{pa 1950 1667}%
\special{pa 1960 1653}%
\special{pa 1965 1647}%
\special{pa 1975 1633}%
\special{pa 1980 1627}%
\special{pa 1990 1613}%
\special{pa 1995 1607}%
\special{pa 2005 1593}%
\special{pa 2010 1587}%
\special{pa 2020 1573}%
\special{pa 2025 1567}%
\special{pa 2035 1553}%
\special{pa 2040 1547}%
\special{pa 2050 1533}%
\special{pa 2055 1527}%
\special{pa 2065 1513}%
\special{pa 2070 1507}%
\special{pa 2080 1493}%
\special{pa 2085 1487}%
\special{pa 2095 1473}%
\special{pa 2100 1467}%
\special{pa 2110 1453}%
\special{pa 2115 1447}%
\special{pa 2125 1433}%
\special{pa 2130 1427}%
\special{pa 2140 1413}%
\special{pa 2145 1407}%
\special{pa 2155 1393}%
\special{pa 2160 1387}%
\special{pa 2170 1373}%
\special{pa 2175 1367}%
\special{pa 2185 1353}%
\special{pa 2190 1347}%
\special{pa 2200 1333}%
\special{pa 2205 1327}%
\special{pa 2215 1313}%
\special{pa 2220 1307}%
\special{pa 2230 1293}%
\special{pa 2235 1287}%
\special{pa 2245 1273}%
\special{pa 2250 1267}%
\special{pa 2260 1253}%
\special{pa 2265 1247}%
\special{pa 2275 1233}%
\special{pa 2280 1227}%
\special{pa 2290 1213}%
\special{pa 2295 1207}%
\special{pa 2300 1200}%
\special{fp}%
% FUNC 2 0 3 0 Black White  
% 10 1400 1200 2600 2400 1400 2400 1600 2400 1400 2200 1400 1200 2600 2400 0 2 0 0 0 0
% 3x/2
\special{pn 8}%
\special{pa 1400 2400}%
\special{pa 1405 2392}%
\special{pa 1415 2378}%
\special{pa 1425 2362}%
\special{pa 1435 2348}%
\special{pa 1445 2332}%
\special{pa 1455 2318}%
\special{pa 1465 2302}%
\special{pa 1475 2288}%
\special{pa 1485 2272}%
\special{pa 1495 2258}%
\special{pa 1505 2242}%
\special{pa 1515 2228}%
\special{pa 1525 2212}%
\special{pa 1535 2198}%
\special{pa 1545 2182}%
\special{pa 1555 2168}%
\special{pa 1565 2152}%
\special{pa 1575 2138}%
\special{pa 1585 2122}%
\special{pa 1595 2108}%
\special{pa 1605 2092}%
\special{pa 1615 2078}%
\special{pa 1625 2062}%
\special{pa 1635 2048}%
\special{pa 1645 2032}%
\special{pa 1655 2018}%
\special{pa 1665 2002}%
\special{pa 1675 1988}%
\special{pa 1685 1972}%
\special{pa 1695 1958}%
\special{pa 1705 1942}%
\special{pa 1715 1928}%
\special{pa 1725 1912}%
\special{pa 1735 1898}%
\special{pa 1745 1882}%
\special{pa 1755 1868}%
\special{pa 1765 1852}%
\special{pa 1775 1838}%
\special{pa 1785 1822}%
\special{pa 1795 1808}%
\special{pa 1805 1792}%
\special{pa 1815 1778}%
\special{pa 1825 1762}%
\special{pa 1835 1748}%
\special{pa 1845 1732}%
\special{pa 1855 1718}%
\special{pa 1865 1702}%
\special{pa 1875 1688}%
\special{pa 1885 1672}%
\special{pa 1895 1658}%
\special{pa 1905 1642}%
\special{pa 1915 1628}%
\special{pa 1925 1612}%
\special{pa 1935 1598}%
\special{pa 1945 1582}%
\special{pa 1955 1568}%
\special{pa 1965 1552}%
\special{pa 1975 1538}%
\special{pa 1980 1530}%
\special{pa 1985 1523}%
\special{pa 1990 1515}%
\special{pa 1995 1508}%
\special{pa 2005 1492}%
\special{pa 2010 1485}%
\special{pa 2015 1477}%
\special{pa 2020 1470}%
\special{pa 2025 1462}%
\special{pa 2035 1448}%
\special{pa 2045 1432}%
\special{pa 2055 1418}%
\special{pa 2065 1402}%
\special{pa 2075 1388}%
\special{pa 2080 1380}%
\special{pa 2085 1373}%
\special{pa 2090 1365}%
\special{pa 2095 1358}%
\special{pa 2105 1342}%
\special{pa 2110 1335}%
\special{pa 2115 1327}%
\special{pa 2120 1320}%
\special{pa 2125 1312}%
\special{pa 2135 1298}%
\special{pa 2145 1282}%
\special{pa 2155 1268}%
\special{pa 2165 1252}%
\special{pa 2175 1238}%
\special{pa 2180 1230}%
\special{pa 2185 1223}%
\special{pa 2190 1215}%
\special{pa 2195 1208}%
\special{pa 2200 1200}%
\special{fp}%
% FUNC 2 0 3 0 Black White  
% 10 1400 1200 2600 2400 1400 2400 1600 2400 1400 2200 1400 1200 2600 2400 0 2 0 0 0 0
% 2x
\special{pn 8}%
\special{pa 1400 2400}%
\special{pa 2000 1200}%
\special{fp}%
% STR 2 0 3 0 Black White  
% 4 1400 700 1400 800 5 0 0 0
% $E_{2}$
\put(14.0000,-8.0000){\makebox(0,0){$E_{2}$}}%
% STR 2 0 3 0 Black White  
% 4 3000 2300 3000 2400 5 0 0 0
% $E_{1}$
\put(30.0000,-24.0000){\makebox(0,0){$E_{1}$}}%
% STR 2 0 3 0 Black White  
% 4 2800 1700 2800 1800 5 0 0 0
% $R$
\put(28.0000,-18.0000){\makebox(0,0){$R$}}%
% STR 2 0 3 0 Black White  
% 4 2800 1500 2800 1600 5 0 0 0
% $P$
\put(28.0000,-16.0000){\makebox(0,0){$P$}}%
% STR 2 0 3 0 Black White  
% 4 2800 900 2800 1000 5 0 0 0
% $S$
\put(28.0000,-10.0000){\makebox(0,0){$S$}}%
% STR 2 0 3 0 Black White  
% 4 2400 900 2400 1000 5 0 0 0
% $Q$
\put(24.0000,-10.0000){\makebox(0,0){$Q$}}%
% STR 2 0 3 0 Black White  
% 4 2200 900 2200 1000 5 0 0 0
% $T$
\put(22.0000,-10.0000){\makebox(0,0){$T$}}%
% STR 2 0 3 0 Black White  
% 4 2000 900 2000 1000 5 0 0 0
% $U$
\put(20.0000,-10.0000){\makebox(0,0){$U$}}%
% LINE 2 0 3 0 Black White  
% 2 1400 1200 2600 2400
% 
\special{pn 8}%
\special{pa 1400 1200}%
\special{pa 2600 2400}%
\special{fp}%
% CIRCLE 2 0 0 0 Black Black  
% 4 4000 1800 4070 1810 4070 1810 4070 1810
% 
\special{sh 1.000}%
\special{ia 4000 1800 71 71 0.0000000 6.2831853}%
\special{pn 8}%
\special{ar 4000 1800 71 71 0.0000000 6.2831853}%
% CIRCLE 2 0 0 0 Black Black  
% 4 6000 1800 6070 1810 6070 1810 6070 1810
% 
\special{sh 1.000}%
\special{ia 6000 1800 71 71 0.0000000 6.2831853}%
\special{pn 8}%
\special{ar 6000 1800 71 71 0.0000000 6.2831853}%
% CIRCLE 2 0 3 0 Black White  
% 4 5000 1800 5070 1810 5070 1810 5070 1810
% 
\special{pn 8}%
\special{ar 5000 1800 71 71 0.0000000 6.2831853}%
% CIRCLE 2 0 3 0 Black White  
% 4 4400 1800 4470 1810 4470 1810 4470 1810
% 
\special{pn 8}%
\special{ar 4400 1800 71 71 0.0000000 6.2831853}%
% CIRCLE 2 0 0 0 Black Black  
% 4 5200 1800 5270 1810 5270 1810 5270 1810
% 
\special{sh 1.000}%
\special{ia 5200 1800 71 71 0.0000000 6.2831853}%
\special{pn 8}%
\special{ar 5200 1800 71 71 0.0000000 6.2831853}%
% CIRCLE 2 0 3 0 Black White  
% 4 5600 1800 5670 1810 5670 1810 5670 1810
% 
\special{pn 8}%
\special{ar 5600 1800 71 71 0.0000000 6.2831853}%
% CIRCLE 2 0 3 0 Black White  
% 4 5380 1800 5450 1810 5450 1810 5450 1810
% 
\special{pn 8}%
\special{ar 5380 1800 71 71 0.0000000 6.2831853}%
% CIRCLE 2 0 0 0 Black Black  
% 4 4720 1800 4790 1810 4790 1810 4790 1810
% 
\special{sh 1.000}%
\special{ia 4720 1800 71 71 0.0000000 6.2831853}%
\special{pn 8}%
\special{ar 4720 1800 71 71 0.0000000 6.2831853}%
% STR 2 0 3 0 Black White  
% 4 4000 1500 4000 1600 5 0 0 0
% $E_{1}$
\put(40.0000,-16.0000){\makebox(0,0){$E_{1}$}}%
% STR 2 0 3 0 Black White  
% 4 6000 1500 6000 1600 5 0 0 0
% $E_{2}$
\put(60.0000,-16.0000){\makebox(0,0){$E_{2}$}}%
% STR 2 0 3 0 Black White  
% 4 5600 1500 5600 1600 5 0 0 0
% $U$
\put(56.0000,-16.0000){\makebox(0,0){$U$}}%
% STR 2 0 3 0 Black White  
% 4 5200 1500 5200 1600 5 0 0 0
% $Q$
\put(52.0000,-16.0000){\makebox(0,0){$Q$}}%
% STR 2 0 3 0 Black White  
% 4 4400 1500 4400 1600 5 0 0 0
% $R$
\put(44.0000,-16.0000){\makebox(0,0){$R$}}%
% STR 2 0 3 0 Black White  
% 4 5000 1500 5000 1600 5 0 0 0
% $S$
\put(50.0000,-16.0000){\makebox(0,0){$S$}}%
% STR 2 0 3 0 Black White  
% 4 4720 1500 4720 1600 5 0 0 0
% $P$
\put(47.2000,-16.0000){\makebox(0,0){$P$}}%
% STR 2 0 3 0 Black White  
% 4 5390 1500 5390 1600 5 0 0 0
% $T$
\put(53.9000,-16.0000){\makebox(0,0){$T$}}%
\end{picture}}%
\end{center}

\caption{A regular simplicial subdivision $\Sigma^*$ which is admissible for $f$.}
\label{exm4_3_regular_subdivision}
\end{figure}

Figure \ref{exm4_3_regular_subdivision} is a regular simplicial subdivision $\Sigma^*$ (of $N^+_\br$) 
which is admissible for $f$, 
where we set 
$$R= \ ^t(2,1),\ S= \ ^t(1,1),\ T= \ ^t(2,3),\ U= \ ^t(1,2).$$
$\Sigma^*$ is called the canonical regular subdivision of $\Gamma^*(f)$ (\cite{Oka1997}, p.137) and 
it is convenient in the sense of Definition \ref{convenient-subdivision-f}. 

All the $2$-dimensional  cones (up to permutations of vertices) of $\Sigma^*$ are as follows:
\begin{equation}\label{7charts}
\begin{array}{l}
\tau_1 := \Cone (E_1,R) = 
\begin{pmatrix}
1 & 2 \\
0 & 1
\end{pmatrix}
,\ \ 
\tau_2 := \Cone (R,P) = 
\begin{pmatrix}
2 & 3 \\
1 & 2
\end{pmatrix}
,\ \ 
\tau_3 := \Cone (P,S) = 
\begin{pmatrix}
3 & 1 \\
2 & 1
\end{pmatrix},\\
\tau_4 := \Cone (S,Q) = 
\begin{pmatrix}
1 & 3 \\
1 & 4
\end{pmatrix}
,\ \ 
\tau_5 := \Cone (Q,T) = 
\begin{pmatrix}
3 & 2 \\
4 & 3
\end{pmatrix}
,\ \ 
\tau_6 := \Cone (T,U) = 
\begin{pmatrix}
2 & 1 \\
3 & 2
\end{pmatrix},\\
\text{and}\ \ 
\displaystyle 
\tau_7 := \Cone (U,E_2) = 
\begin{pmatrix}
1 & 0 \\
2 & 1
\end{pmatrix}.
\end{array}
\end{equation}

We consider the toric modification 
\begin{equation}\label{modification_exm4_3}
\hat{\pi} : X \to \bc^2
\end{equation}
associated with $\Sigma^*$, 
which is called the canonical toric modification for $f$ (\cite{Oka1997}, p.138). 

The face functions of the germ $(f,\zero)$ are
$$
f_U = f_T = z_1^4\bar{z}_1^2,\ \ f_S = \bar{z}_1^2z_2^3,\ \ f_R = z_2^3\bar{z}_2^3,\ \ 
f_Q = z_1^4\bar{z}_1^2 + \bar{z}_1^2z_2^3,\ \ f_P = \bar{z}_1^2z_2^3 + z_2^3\bar{z}_2^3 .
$$
The radial and polar degrees of these face functions are calculated in Table \ref{table4323}. 

\begin{table}
\begin{center}
%\begin{tabular}{c}
%\begin{minipage}{  0.3\hsize}
%\begin{center}
\begin{tabular}{|c|c|c|}
\hline
$0$-dim face $\De (U)$       & \multicolumn{2}{|c|}{$f_{U}=z_{1}^{4}\overline{z_{1}}^{2}$}   \\ \hline
\multirow{2}{*}{weight vector:} & $1$     &  $4+2, 4-2$   \\ \cline{3-3} 
                                           & $2$     &  $0+0, 0-0$  \\ \hline
radial degree                   &         & $6$ \\ \hline
polar degree                    &         & $2$ \\ \hline
\end{tabular}
%\end{center}
%\end{minipage}
\ \ \ 
%\begin{minipage}{  0.3\hsize}
%\begin{center}
\begin{tabular}{|c|c|c|}
\hline
$0$-dim face $\De (T)$              & \multicolumn{2}{|c|}{$f_{T}=z_{1}^{4}\overline{z_{1}}^{2}$} \\ \hline
\multirow{2}{*}{weight vector:} & $2$     &   $4+2, 4-2$   \\ \cline{3-3} 
                                           & $3$     &   $0+0, 0-0$    \\ \hline
radial degree                   &         & $12$ \\ \hline
polar degree                    &         & $4$ \\ \hline
\end{tabular}
%\end{center}
%\end{minipage}
%\end{tabular}
\end{center}

\medskip

\begin{center}
%\begin{tabular}{c}
%\begin{minipage}{  0.3\hsize}
%\begin{center}
\begin{tabular}{|c|c|c|}
\hline
$0$-dim face $\De (S)$     & \multicolumn{2}{|c|}{$f_{S}=z_{2}^{3}\overline{z_{1}}^{2}$} \\ \hline
\multirow{2}{*}{weight vector:} & $1$     & $0+2, 0-2$      \\ \cline{3-3} 
                                           & $1$     & $3+0, 3-0$      \\ \hline
radial degree                   &         & $5$ \\ \hline
polar degree                    &         & $1$ \\ \hline
\end{tabular}
%\end{center}
%\end{minipage}
\ \ \ 
%\begin{minipage}{  0.3\hsize}
%\begin{center}
\begin{tabular}{|c|c|c|}
\hline
$0$-dim face $\De (R)$           & \multicolumn{2}{|c|}{$f_{R}=z_{2}^{3}\overline{z_{2}}^{3}$} \\ \hline
\multirow{2}{*}{weight vector:} & $2$     &  $0+0,0-0$        \\ \cline{3-3} 
                                           & $1$     &  $3+3,3-3$  \\ \hline
radial degree                   &         & $6$ \\ \hline
polar degree                    &         & $0$ \\ \hline
\end{tabular}
%\end{center}
%\end{minipage}
%\end{tabular}
\end{center}

\medskip

\begin{center}
\begin{tabular}{|c|c|c|c|}
\hline
$1$-dim face $\De (Q)$         & \multicolumn{3}{|c|}{$f_{Q}=z_{1}^{4}\overline{z_{1}}^{2} + z_{2}^{3}\overline{z_{1}}^{2}$} \\ \hline
\multirow{2}{*}{weight vector:} & $3$     &  $4+2,4-2$       & $0+2,0-2$               \\ \cline{3-4} 
                                           & $4$     &   $0+0,0-0$      & $3+0,3-0$               \\ \hline
                                           &           & the first term    & the second term \\ \hline
radial degree                   &                      & $18$           & $18$            \\ \hline
polar degree                    &                      & $6$            & $6$             \\ \hline
\end{tabular}
\end{center}

\medskip

\begin{center}
\begin{tabular}{|c|c|c|c|}
\hline
$1$-dim face $\De (P)$          & \multicolumn{3}{|c|}{$f_{P}=z_{2}^{3}\overline{z_{1}}^{2} + z_{2}^{3}\overline{z_{2}}^{3}$} \\ \hline
\multirow{2}{*}{weight vector:} & $3$     &   $0+2,0-2$      &    $0+0,0-0$            \\ \cline{3-4} 
                                           & $2$     &   $3+0,3-0$      &     $3+3,3-3$           \\ \hline
                                              &         & the first term & the second term \\ \hline
radial degree                   &                  & $12$             & $12$            \\ \hline
polar degree                    &                  & $0$               & $0$             \\ \hline
\end{tabular}
\end{center}

\bigskip

\caption{The face functions of the germ $((z_1^4 + z_2^3)\overline{(z_1^2 + z_2^3)},\zero)$.}
\label{table4323}
\end{table}

Let us discuss the Newton non-degeneracy of the germ $(f,\zero)$. 
For weight vectors $P= {}^t(3,2), \ Q= {}^t(3,4)$, $\dim \De(P) = \dim \De(Q) = 1$, and 
$$
f_{P}(\bbz,\overline{\bbz})=z_{2}^{3}\overline{z_{1}}^{2}+z_{2}^{3}\overline{z_{2}}^{3}, \ \ \ \ 
f_{Q}(\bbz,\overline{\bbz})=z_{1}^{4}\overline{z_{1}}^{2}+z_{2}^{3}\overline{z_{1}}^{2} .
$$

For $P$, we have 
$$
\frac{\partial f_{P}}{\partial z_{1}}=0, \ \ \ \frac{\partial f_{P}}{\partial \overline{z_{1}}}=2z_{2}^{3}\overline{z_{1}} .
$$
Hence, we have 
$$
\left|\frac{\partial f_{P}}{\partial z_{1}}\right| \neq \left|\frac{\partial f_{P}}{\partial \overline{z_{1}}}\right|
$$
on $\bc^{*2}$. 
Then by Corollary \ref{regular-criterion}, 
%%%Proposition \ref{critical}, 
$f_{P}$ has no mixed critical point on $\bc^{*2}$. 

The function $f_{P}(z_{1},1)=\overline{z_{1}}^{2}+1$ is a surjective map from $\bc^*$ onto $\bc \setminus \{1\}$. 
Since $f_{P}(\sqrt{\frac{3}{2}},\sqrt[3]{\frac{1}{2}})=1$, we see $f_{P}:\bc^{*2} \to \bc$ is surjective. 
Thus $(f,\zero)$ is strongly Newton non-degenerate over $\De(P)$. 

For $Q$, we have 
$$
\frac{\partial f_{Q}}{\partial z_{2}}=3z_{2}^{2}\overline{z_{1}}^{2}, \ \ \ \frac{\partial f_{Q}}{\partial \overline{z_{2}}}=0 .
$$
Hence, we also have 
$$
\left|\frac{\partial f_{Q}}{\partial z_{2}}\right| \neq \left|\frac{\partial f_{Q}}{\partial \overline{z_{2}}}\right|
%%%% \left|\frac{\partial f_{Q}}{\partial z_{1}}\right|-\left|\frac{\partial f_{Q}}{\partial \overline{z_{1}}}\right|=-4|z_{1}|^{2}|z_{2}|^{6}
$$
on $\bc^{*2}$. 
Then by Corollary \ref{regular-criterion}, 
$f_{Q}$ has no mixed critical point on $\bc^{*2}$. 

The function $f_{Q}(1,z_{2})=1+z_{2}^{3}$ is a surjective map from $\bc^*$ onto $\bc \setminus \{1\}$. 
Since $f_{Q}(\sqrt{\frac{1}{2}},\sqrt[3]{\frac{7}{4}})=1$, we see $f_{Q}:\bc^{*2} \to \bc$ is surjective. 
Thus $(f,\zero)$ is strongly Newton non-degenerate over $\De(Q)$. 

For $0$-dimensional faces 
$\De (U) = \De (T) = \{(6,0)\},\ \De (S) = \{(2,3)\}$ and $\De (R) = \{(0,6)\}$, 
the corresponding face functions are 
$$f_U = f_T = z_{1}^{4}\overline{z_{1}}^{2}, \ \ \ f_S = z_{2}^{3}\overline{z_{1}}^{2}, \ \ \ f_R = z_{2}^{3}\overline{z_{2}}^{3}.$$
Since 
$$
\begin{array}{ll}
\frac{\partial}{\partial z_{1}}z_{1}^{4}\overline{z_{1}}^{2} = 4z_{1}^{3}\overline{z_{1}}^{2}, & \frac{\partial}{\partial \overline{z_{1}}}z_{1}^{4}\overline{z_{1}}^{2} = 2z_{1}^{4}\overline{z_{1}}, \\
\frac{\partial}{\partial z_{1}}z_{2}^{3}\overline{z_{1}}^{2} = 0, \ \ \ \text{and}\            & \frac{\partial}{\partial \overline{z_{1}}}z_{2}^{3}\overline{z_{1}}^{2} = 2z_{2}^{3}\overline{z_{1}}, 
%% \frac{\partial}{\partial z_{2}}z_{2}^{3}\overline{z_{2}}^{3} = 3z_{2}^{2}\overline{z_{2}}^{3}, & \frac{\partial}{\partial \overline{z_{2}}}z_{2}^{3}\overline{z_{2}}^{3} = 3z_{2}^{3}\overline{z_{2}}^{2} 
\end{array}
$$
$(f,\zero)$ is strongly Newton non-degenerate over $\De(U)=\De(T)$ and $\De(S)$. 

However, the face function $z_{2}^{3}\overline{z_{2}}^{3} = |z_2|^6$ takes only real values. 
%%$$
%%\left|\frac{\partial}{\partial z_{2}}z_{2}^{3}\overline{z_{2}}^{3}\right|-\left|\frac{\partial}{\partial \overline{z_{2}}}z_{2}^{3}\overline{z_{2}}^{3}\right|=0
%%$$
Hence, every point $(z_1, z_2) \in \bc^2$ is a mixed critical point of $z_{2}^{3}\overline{z_{2}}^{3}$. 
Since $z_{2}^{3}\overline{z_{2}}^{3} = |z_2|^6 >0$ on $\bc^{*2}$, $0$ is not a mixed critical value of 
$z_{2}^{3}\overline{z_{2}}^{3} : \bc^{*2} \to \bc$. 
Thus $(f,\zero)$ is not strongly Newton non-degenerate but Newton non-degenerate over $\De(R)$. 

The second term $z_1^4\bar{z}_2^3$ of the mixed polynomial \eqref{exam4-3-mixed} belongs to 
the interior of the (radial) Newton polyhedron of $f(\bbz, \bar{\bbz})$. 
Let us replace $z_1^4\bar{z}_2^3$ by $z_1^a \bar{z}_1^{4-a} z_2^b \bar{z}_2^{3-b}$, and 
consider the mixed polynomials 
$$
f_{a,b}(\bbz, \bar{\bbz}) := z_1^4 \bar{z}_1^2 + z_1^a \bar{z}_1^{4-a} z_2^b \bar{z}_2^{3-b} + \bar{z}_1^2 z_2^3 + z_2^3 \bar{z}_2^3,
$$
where $0 \leq a \leq 4, \ 0\leq b \leq 3$ are integers. 
Then, $(f_{a,b},\zero)$ is convenient. %%

Every mixed polynomial germ $(f_{a,b},\zero)$ has the same (radial) Newton polyhedron as $(f,\zero)$. 
Moreover, it has the same face functions as $(f,\zero)$. %%
Hence, the germ 
$(f_{a,b},\zero)$ is also a Newton non-degenerate mixed polynomial germ %%
of strongly polar non-negative mixed weighted homogeneous face type. 

\medskip

We put $f := f_{a,b}$ hereafter. 

\medskip

%% k=1 %%%%%%%%%%%%%%%%%%%%%%%%%%%%%%%%%%%%%%%%%%%%%%%%%%%%%%%%%%%%
Take a toric chart 
$\sigma = \Cone (P_{1},P_{2})$, where $P_{j}= \ ^t(p_{1j},p_{2j}) \ (j=1,2)$ and 
$P_{1}$ is strictly positive (namely, $P_{1} \neq E_{1},E_{2}$). %%

Then, for $(u_{1},u_{2}) \in \bc_\sigma^2$, %%
$$
\pi_{\sigma}(u_{1},u_{2})=(u_{1}^{p_{11}}u_{2}^{p_{12}},u_{1}^{p_{21}}u_{2}^{p_{22}})
$$
and we have 
$$
\begin{array}{cl}
   &  \pi_{\sigma}^{*}f(u_{1},u_{2})\\
= & u_{1}^{4p_{11}} \overline{u_{1}}^{2p_{11}} u_{2}^{4p_{12}} \overline{u_{2}}^{2p_{12}}\, + \,u_{1}^{ap_{11}+bp_{21}} \overline{u_{1}}^{(4-a)p_{11}+(3-b)p_{21}} u_{2}^{ap_{12}+bp_{22}} \overline{u_{2}}^{(4-a)p_{12}+(3-b)p_{22}} \\
   & + \, u_{1}^{3p_{21}} \overline{u_{1}}^{2p_{11}} u_{2}^{3p_{22}} \overline{u_{2}}^{2p_{12}}\, +\, u_{1}^{3p_{21}} \overline{u_{1}}^{3p_{21}} u_{2}^{3p_{22}} \overline{u_{2}}^{3p_{22}} .
\end{array}
$$

We set $\tau := \Cone (P_{1})$. 

We have 
$$
f_{P_1}(z_{1},z_{2}) = c_1z_{1}^{4}\overline{z_{1}}^{2} + c_2\overline{z_{1}}^{2}z_{2}^{3} + c_3z_{2}^{3}\overline{z_{2}}^{3} ,
$$
where 
$$
(c_1, c_2, c_3) = 
\begin{cases}
(1,0,0) & (\text{if}\ P_1 = U \ \text{or} \ T)\\
(1,1,0) & (\text{if}\ P_1 = Q)\\
(0,1,0) & (\text{if}\ P_1 = S)\\
(0,1,1) & (\text{if}\ P_1 = P)\\
(0,0,1) & (\text{if}\ P_1 = R) .
\end{cases}
$$

Hence, we have 
$$
\rdeg f_{P_{1}}=\begin{cases}
6p_{11} & (c_{1}=1) \\
2p_{11}+3p_{21} & (c_{2}=1) \\
6p_{21} & (c_{3}=1) ,
\end{cases}
\ \ \ \ \ 
\pdeg f_{P_{1}}=\begin{cases}
2p_{11} & (c_{1}=1) \\
-2p_{11}+3p_{21} & (c_{2}=1) \\
0 & (c_{3}=1) ,
\end{cases}
$$
and 
\begin{eqnarray*}
\frac{\rdeg f_{P_{1}}+\pdeg f_{P_{1}}}{2} &=& 
\begin{cases}
4p_{11} & (c_{1}=1) \\
3p_{21} & (c_{2}=1 \ \text{or} \ c_{3}=1)
\end{cases} \\ 
\frac{\rdeg f_{P_{1}}-\pdeg f_{P_{1}}}{2} &=& 
\begin{cases}
2p_{11} & (c_{1}=1 \ \text{or} \ c_{2}=1) \\
3p_{21} & (c_{3}=1) .
\end{cases}
\end{eqnarray*}

The pull-back of the face function $f_{P_{1}}$ by $\pi_{\sigma}$ is as follows:
\begin{eqnarray*}
\pi_{\sigma}^{*}f_{P_{1}}(u_{1},u_{2}) = c_{1}u_{1}^{4p_{11}}\overline{u_{1}}^{2p_{11}}u_{2}^{4p_{12}}\overline{u_{2}}^{2p_{12}} + c_{2}u_{1}^{3p_{21}}\overline{u_{1}}^{2p_{11}}u_{2}^{3p_{22}}\overline{u_{2}}^{2p_{12}}+c_{3}u_{1}^{3p_{21}}\overline{u_{1}}^{3p_{21}}u_{2}^{3p_{22}}\overline{u_{2}}^{3p_{22}}
\end{eqnarray*}

\bigskip

\noindent
Case\ $c_{1}=1$: %%
We have
\begin{eqnarray*}
\pi_{\sigma}^{*}f(u_{1},u_{2}) &=& 
u_{1}^{4p_{11}}\overline{u_{1}}^{2p_{11}} %%
\left( u_{1}^{(3p_{21}-4p_{11})}u_{2}^{3p_{22}}\overline{u_{2}}^{2p_{12}}  \right .
 + u_{1}^{(3p_{21}-4p_{11})} \overline{u_{1}}^{(3p_{21}-2p_{11})} u_{2}^{3p_{22}} \overline{u_{2}}^{3p_{22}} \\
& &  + u_{2}^{4p_{12}}\overline{u_{2}}^{2p_{12}} 
+ \left . \underline{u_{1}^{((a-4)p_{11}+bp_{21})} \overline{u_{1}}^{((2-a)p_{11}+(3-b)p_{21})} u_{2}^{(ap_{12}+bp_{22})} \overline{u_{2}}^{((4-a)p_{12}+(3-b)p_{22})} } \right) .
\end{eqnarray*}
Here 
we have %%
$6p_{11} = 2p_{11}+3p_{21}$  %%
or 
$6p_{11} < 2p_{11}+3p_{21}$.  %%
Hence, 
$3p_{21}-2p_{11}>3p_{21}-4p_{11} \geq 0$. 

\bigskip

\noindent
Case\ $c_{2}=1$: %%
We have
\begin{eqnarray*}
\pi_{\sigma}^{*}f(u_{1},u_{2}) &=& 
u_{1}^{3p_{21}}\overline{u_{1}}^{2p_{11}} %%
\left( u_{1}^{(4p_{11}-3p_{21})}u_{2}^{4p_{12}}\overline{u_{2}}^{2p_{12}} \right .
+ \overline{u_{1}}^{(3p_{21}-2p_{11})} u_{2}^{3p_{22}} \overline{u_{2}}^{3p_{22}} \\
 & & + u_{2}^{3p_{22}}\overline{u_{2}}^{2p_{12}}
+ \left . \underline{u_{1}^{(ap_{11}+(b-3)p_{21})} \overline{u_{1}}^{((2-a)p_{11}+(3-b)p_{21})} u_{2}^{(ap_{12}+bp_{22})} \overline{u_{2}}^{((4-a)p_{12}+(3-b)p_{22})} } \right) .
\end{eqnarray*}
Here 
we have %%
$2p_{11}+3p_{21} \leq 6p_{11}$ and $2p_{11}+3p_{21} \leq 6p_{21}$. 
Hence, 
$4p_{11}-3p_{21} \geq 0$ and $3p_{21}-2p_{11} \geq 0$. 

\bigskip

\noindent
Case\ $c_{3}=1$:  %%
We have
\begin{eqnarray*}
\pi_{\sigma}^{*}f(u_{1},u_{2}) &=&
u_{1}^{3p_{21}}\overline{u_{1}}^{3p_{21}}  %%
\left( \overline{u_{1}}^{(2p_{11}-3p_{21})} u_{2}^{3p_{22}} \overline{u_{2}}^{2p_{12}} \right .
+ u_{1}^{(4p_{11}-3p_{21})} \overline{u_{1}}^{(2p_{11}-3p_{21})} u_{2}^{3p_{22}} \overline{u_{2}}^{3p_{22}} \\
& & + u_{2}^{3p_{22}}\overline{u_{2}}^{3p_{22}}
+ \left . \underline{u_{1}^{(ap_{11}+(b-3)p_{21})} \overline{u_{1}}^{((4-a)p_{11}-bp_{21})} u_{2}^{(ap_{12}+bp_{22})}\overline{u_{2}}^{((4-a)p_{12}+(3-b)p_{22})} } \right) .
\end{eqnarray*}
Here 
we have %%
$6p_{21} \leq 2p_{11}+3p_{21}$. 
Hence, 
$2p_{11}-3p_{21} \geq 0$ and $4p_{11}-3p_{21} \geq 2p_{11} > 0$. 

For the exponents of the underlined terms above, we calculate 
\begin{eqnarray*}
c_{1}=1 &\Longrightarrow& ((a-4)p_{11}+bp_{21}) + ((2-a)p_{11}+(3-b)p_{21}) = -2p_{11}+3p_{21} \geq 2p_{11} \geq 2 \\
c_{2}=1 &\Longrightarrow& (ap_{11}+(b-3)p_{21}) + ((2-a)p_{11}+(3-b)p_{21}) = 2p_{11} \geq 2 \\
c_{3}=1 &\Longrightarrow& (ap_{11}+(b-3)p_{21}) + ((4-a)p_{11}-bp_{21}) = 4p_{11}-3p_{21} \geq 2p_{11} \geq 2 . 
\end{eqnarray*}

Hence, we may define the value of each underlined term at $(0, u_2)$ to be $0$.  %%
Then each underlined term is of class $C^{1}$. 
Hence $\tilde{R}, \tilde{f}$ are also of class $C^{1}$. 

\medskip

Let us consider the case $P_1 = R = {}^t(2,1)$. Namely, $(c_1, c_2, c_3) = (0, 0, 1)$ and $(p_{11},p_{21})=(2,1)$. 
Recall that $f_R = z_{2}^{3}\overline{z_{2}}^{3}$. 
We have 
$$
\rdeg f_{P_{1}}= 6
\ \ \ \text{and}\ \ \ 
\pdeg f_{P_{1}}= 0 ,
$$
and hence, 
$$
\frac{\rdeg f_{P_{1}}+\pdeg f_{P_{1}}}{2} = 3
\ \ \ \text{and}\ \ \ 
\frac{\rdeg f_{P_{1}}-\pdeg f_{P_{1}}}{2} = 3 .
$$

We have
\begin{eqnarray*}
\pi_{\sigma}^{*}f(u_{1},u_{2}) &=&
u_{1}^{3}\overline{u_{1}}^{3}  %%
\left( \overline{u_{1}}^{(4-3)} u_{2}^{3p_{22}} \overline{u_{2}}^{2p_{12}} \right .
+ u_{1}^{(8-3)} \overline{u_{1}}^{(4-3)} u_{2}^{3p_{22}} \overline{u_{2}}^{3p_{22}} \\
& & + u_{2}^{3p_{22}}\overline{u_{2}}^{3p_{22}}
+ \left . \underline{u_{1}^{(2a+(b-3))} \overline{u_{1}}^{(2(4-a)-b)} u_{2}^{(ap_{12}+bp_{22})}\overline{u_{2}}^{((4-a)p_{12}+(3-b)p_{22})} } \right) \\
 &=&
u_{1}^{3}\overline{u_{1}}^{3}  %%
\left( \overline{u_{1}} u_{2}^{3p_{22}} \overline{u_{2}}^{2p_{12}} \right .
+ u_{1}^5 \overline{u_{1}} u_{2}^{3p_{22}} \overline{u_{2}}^{3p_{22}} \\
& & + u_{2}^{3p_{22}}\overline{u_{2}}^{3p_{22}}
+ \left . \underline{u_{1}^{(2a+(b-3))} \overline{u_{1}}^{(2(4-a)-b)} u_{2}^{(ap_{12}+bp_{22})}\overline{u_{2}}^{((4-a)p_{12}+(3-b)p_{22})} } \right) 
\end{eqnarray*}
$$
= 
\begin{cases}
u_{1}^{3}\overline{u_{1}}^{3} \left( \overline{u_{1}} \overline{u_{2}}^2 + u_{1}^5 \overline{u_{1}} + 1 + \underline{u_{1}^{(2a+(b-3))} \overline{u_{1}}^{(2(4-a)-b)} u_{2}^a \overline{u_{2}}^{(4-a)} } \right) 
\ \ \ \text{if} \ P_2 = E_1 = {}^t(1,0) ,\\
u_{1}^{3}\overline{u_{1}}^{3} \left( \overline{u_{1}} u_{2}^6 \overline{u_{2}}^6 + u_{1}^5 \overline{u_{1}} u_{2}^6 \overline{u_{2}}^6 + u_{2}^6 \overline{u_{2}}^6 + \underline{u_{1}^{(2a+(b-3))} \overline{u_{1}}^{(2(4-a)-b)} u_{2}^{(3a+2b)}\overline{u_{2}}^{(3(4-a)+2(3-b))} } \right) \\
\hspace{10.2cm} \ \ \ \text{if} \ P_2 = P   ={}^t(3,2) .
\end{cases}
$$

Hence, we have 
$$
\begin{cases}
\tilde{f}_R = 1, \ \tilde{R} = \overline{u_{1}} \overline{u_{2}}^2 + u_{1}^5 \overline{u_{1}} + \underline{u_{1}^{(2a+(b-3))} \overline{u_{1}}^{(2(4-a)-b)} u_{2}^a \overline{u_{2}}^{(4-a)} } 
\ \ \ \text{if} \ P_2 = E_1 = {}^t(1,0) ,\\
\tilde{f}_R = u_{2}^6 \overline{u_{2}}^6, \ \tilde{R} = \overline{u_{1}} u_{2}^6 \overline{u_{2}}^6 + u_{1}^5 \overline{u_{1}} u_{2}^6 \overline{u_{2}}^6 + \underline{u_{1}^{(2a+(b-3))} \overline{u_{1}}^{(2(4-a)-b)} u_{2}^{(3a+2b)}\overline{u_{2}}^{(3(4-a)+2(3-b))} } \\
\hspace{10.2cm} \ \ \ \text{if} \ P_2 = P   ={}^t(3,2) .
\end{cases}
$$

Thus we have $\tilde{f}(0, u_2) \neq 0$ if $u_2 \neq 0$. 
Hence, we see that $\tilde V(R)^* = \emptyset$. 

Moreover, in the toric chart $\sigma = \Cone (R,E_1)$, 
we also see that $\tilde V = \emptyset$ in some neighborhood of $\hat{E}(R)$ 
since $\tilde{f}(0, u_2) = \tilde{f}_R(0, u_2) = 1$. %%%%

\medskip

Next let us consider the case $P_1 = U = {}^t(1,2)$. 
Namely, $(c_1, c_2, c_3) = (1, 0, 0)$ and $(p_{11},p_{21})=(1,2)$. 
Recall that $f_U = z_{1}^{4}\overline{z_{1}}^{2}$. 
We have 
$$
\rdeg f_{P_{1}}= 6\ \ \text{and}\ \ 
\pdeg f_{P_{1}}= 2 ,
$$
and hence, 
$$
\frac{\rdeg f_{P_{1}}+\pdeg f_{P_{1}}}{2} = 4\ \ \text{and}\ \ 
\frac{\rdeg f_{P_{1}}-\pdeg f_{P_{1}}}{2} = 2 .
$$

We have
\begin{eqnarray*}
\pi_{\sigma}^{*}f(u_{1},u_{2}) &=& 
u_{1}^{4}\overline{u_{1}}^{2} %%
\left( u_{1}^{(6-4)}u_{2}^{3p_{22}}\overline{u_{2}}^{2p_{12}}  \right .
 + u_{1}^{(6-4)} \overline{u_{1}}^{(6-2)} u_{2}^{3p_{22}} \overline{u_{2}}^{3p_{22}} \\
& &  + u_{2}^{4p_{12}}\overline{u_{2}}^{2p_{12}} 
+ \left . \underline{u_{1}^{((a-4) + 2b)} \overline{u_{1}}^{((2-a)+(3-b)\times 2)} u_{2}^{(ap_{12}+bp_{22})} \overline{u_{2}}^{((4-a)p_{12}+(3-b)p_{22})} } \right)
\end{eqnarray*}
$$
= 
\begin{cases}
u_{1}^{4}\overline{u_{1}}^{2} %%
\left( u_{1}^2u_{2}^{3}
+ u_{1}^2 \overline{u_{1}}^4 u_{2}^3 \overline{u_{2}}^3
+ 1 + \underline{u_{1}^{((a-4) + 2b)} \overline{u_{1}}^{((2-a)+(3-b)\times 2)} u_{2}^{b} \overline{u_{2}}^{(3-b)} } \right)
\ \ \ \text{if} \ P_2 = E_2 = {}^t(0,1) ,\\
u_{1}^{4}\overline{u_{1}}^{2} %%
\left( u_{1}^2u_{2}^{9}\overline{u_{2}}^{4}
+ u_{1}^{2} \overline{u_{1}}^{4} u_{2}^{9} \overline{u_{2}}^{9}
+ u_{2}^{8}\overline{u_{2}}^{4} 
+ \underline{u_{1}^{((a-4) + 2b)} \overline{u_{1}}^{((2-a)+(3-b)\times 2)} u_{2}^{(2a + 3b)} \overline{u_{2}}^{((4-a)\times 2 + (3-b)\times 3)} } \right) \\
\hspace{12.2cm} \ \ \ \text{if} \ P_2 = T ={}^t(2,3) .
\end{cases}
$$

Hence, we have 
$$
\begin{cases}
\tilde{f}_U = 1, \ \tilde{R} = u_{1}^2u_{2}^{3}
+ u_{1}^2 \overline{u_{1}}^4 u_{2}^3 \overline{u_{2}}^3
+ \underline{u_{1}^{((a-4) + 2b)} \overline{u_{1}}^{((2-a)+(3-b)\times 2)} u_{2}^{b} \overline{u_{2}}^{(3-b)} } 
\ \ \ \text{if} \ P_2 = E_2 = {}^t(0,1) ,\\
\tilde{f}_U = u_{2}^8 \overline{u_{2}}^4, \ \tilde{R} = u_{1}^2u_{2}^{9}\overline{u_{2}}^{4}
+ u_{1}^{2} \overline{u_{1}}^{4} u_{2}^{9} \overline{u_{2}}^{9}
+ \underline{u_{1}^{((a-4) + 2b)} \overline{u_{1}}^{((2-a)+(3-b)\times 2)} u_{2}^{(2a + 3b)} \overline{u_{2}}^{((4-a)\times 2 + (3-b)\times 3)} } \\
\hspace{12.2cm} \ \ \ \text{if} \ P_2 = T ={}^t(2,3) .
\end{cases}
$$

Thus we have $\tilde{f}(0, u_2) \neq 0$ if $u_2 \neq 0$. 
Hence, we see that $\tilde V(U)^* = \emptyset$. 

Moreover, in the toric chart $\sigma = \Cone (U,E_2)$, 
we also see that $\tilde V = \emptyset$ in some neighborhood of $\hat{E}(U)$ 
since $\tilde{f}(0, u_2) = \tilde{f}_U(0, u_2) = 1$. %%%

\bigskip

%% k=2 %%%%%%%%%%%%%%%%%%%%%%%%%%%%%%%%%%%%%%%%%%%%%%%%%%%%%%%%%%%%%%%%%%%%%%%%%%%%%%
Now take a toric chart $\tau := \Cone (P_{1}, P_{2})$, 
where $P_{j}= \ ^t(p_{1j},p_{2j}) \ (j=1,2)$, and 
both $P_{1}$ and $P_{2}$ are strictly positive weight vectors (namely, $P_{j} \neq E_{1},E_{2}  \ (j=1,2)$). 

We define the toric chart 
$\sigma := \tau$. %%% k=2

Then there are five such $2$-dimentional charts $\sigma$ as follows (recall the notation \eqref{7charts}):
%%%%
$$
\begin{array}{ll}
\tau_{2} = \Cone{(R,P)}, \ \ & \pi_{\tau_{2}}(u_{1},u_{2})=(u_{1}^{2}u_{2}^{3},u_{1}u_{2}^{2}), \\
\tau_{3} = \Cone{(P,S)}, \ \ & \pi_{\tau_{3}}(u_{1},u_{2})=(u_{1}^{3}u_{2},u_{1}^{2}u_{2}), \\
\tau_{4} = \Cone{(S,Q)}, \ \ & \pi_{\tau_{4}}(u_{1},u_{2})=(u_{1}u_{2}^{3},u_{1}u_{2}^{4}), \\
\tau_{5} = \Cone{(Q,T)}, \ \ & \pi_{\tau_{5}}(u_{1},u_{2})=(u_{1}^{3}u_{2}^{2},u_{1}^{4}u_{2}^{3}), \\
\tau_{6} = \Cone{(T,U)}, \ \ & \pi_{\tau_{6}}(u_{1},u_{2})=(u_{1}^{2}u_{2},u_{1}^{3}u_{2}^{2}) .
\end{array}
$$

Hence, we have
\begin{eqnarray*}
\pi_{\tau_{2}}^{*}f(u_{1},u_{2}) &=& u_{1}^{8}u_{2}^{12}\overline{u_{1}}^{4}\overline{u_{2}}^{6}+u_{1}^{2a+b}u_{2}^{3a+2b}\overline{u_{1}}^{11-2a-b}\overline{u_{2}}^{18-3a-2b}+u_{1}^{3}u_{2}^{6}\overline{u_{1}}^{4}\overline{u_{2}}^{6}+u_{1}^{3}u_{2}^{6}\overline{u_{1}}^{3}\overline{u_{2}}^{6}, \\
\pi_{\tau_{3}}^{*}f(u_{1},u_{2}) &=& u_{1}^{12}u_{2}^{4}\overline{u_{1}}^{6}\overline{u_{2}}^{2}+u_{1}^{3a+2b}u_{2}^{a+b}\overline{u_{1}}^{18-3a-2b}\overline{u_{2}}^{7-a-b}+u_{1}^{6}u_{2}^{3}\overline{u_{1}}^{6}\overline{u_{2}}^{2}+u_{1}^{6}u_{2}^{3}\overline{u_{1}}^{6}\overline{u_{2}}^{3}, \\
\pi_{\tau_{4}}^{*}f(u_{1},u_{2}) &=& u_{1}^{4}u_{2}^{12}\overline{u_{1}}^{2}\overline{u_{2}}^{6}+u_{1}^{a+b}u_{2}^{3a+4b}\overline{u_{1}}^{7-a-b}\overline{u_{2}}^{24-3a-4b}+u_{1}^{3}u_{2}^{12}\overline{u_{1}}^{2}\overline{u_{2}}^{6}+u_{1}^{3}u_{2}^{12}\overline{u_{1}}^{3}\overline{u_{2}}^{12}, \\
\pi_{\tau_{5}}^{*}f(u_{1},u_{2}) &=& u_{1}^{12}u_{2}^{8}\overline{u_{1}}^{6}\overline{u_{2}}^{4}+u_{1}^{3a+2b}u_{2}^{2a+3b}\overline{u_{1}}^{24-3a-4b}\overline{u_{2}}^{17-2a-3b}+u_{1}^{12}u_{2}^{9}\overline{u_{1}}^{6}\overline{u_{2}}^{4}+u_{1}^{12}u_{2}^{9}\overline{u_{1}}^{12}\overline{u_{2}}^{9}, \\
\pi_{\tau_{6}}^{*}f(u_{1},u_{2}) &=& u_{1}^{8}u_{2}^{4}\overline{u_{1}}^{4}\overline{u_{2}}^{2}+u_{1}^{2a+3b}u_{2}^{a+2b}\overline{u_{1}}^{17-2a-3b}\overline{u_{2}}^{10-a-2b}+u_{1}^{9}u_{2}^{6}\overline{u_{1}}^{4}\overline{u_{2}}^{2}+u_{1}^{9}u_{2}^{6}\overline{u_{1}}^{9}\overline{u_{2}}^{6} ,
\end{eqnarray*}
and
%%%%%
\begin{eqnarray*}
\pi_{\tau_{2}}^{*}f(u_{1},u_{2}) &=& u_{1}^{3}u_{2}^{6}\overline{u_{1}}^{3}\overline{u_{2}}^{6}(u_{1}^{5}u_{2}^{6}\overline{u_{1}}+\underline{u_{1}^{2a+b-3}u_{2}^{3a+2b-6}\overline{u_{1}}^{8-2a-b}\overline{u_{2}}^{12-3a-2b}}+\overline{u_{1}}+1), \\
\pi_{\tau_{3}}^{*}f(u_{1},u_{2}) &=& u_{1}^{6}u_{2}^{3}\overline{u_{1}}^{6}\overline{u_{2}}^{2}(u_{1}^{6}u_{2}+\underline{u_{1}^{3a+2b-6}u_{2}^{a+b-3}\overline{u_{1}}^{12-3a-2b}\overline{u_{2}}^{5-a-b}}+1+\overline{u_{2}}), \\
\pi_{\tau_{4}}^{*}f(u_{1},u_{2}) &=& u_{1}^{3}u_{2}^{12}\overline{u_{1}}^{2}\overline{u_{2}}^{6}(u_{1}+\underline{u_{1}^{a+b-3}u_{2}^{3a+4b-12}\overline{u_{1}}^{5-a-b}\overline{u_{2}}^{18-3a-4b}}+1+\overline{u_{1}}\overline{u_{2}}^{6}), \\
\pi_{\tau_{5}}^{*}f(u_{1},u_{2}) &=& u_{1}^{12}u_{2}^{8}\overline{u_{1}}^{6}\overline{u_{2}}^{4}(1+\underline{u_{1}^{3a+2b-12}u_{2}^{2a+3b-8}\overline{u_{1}}^{18-3a-4b}\overline{u_{2}}^{13-2a-3b}}+u_{2}+u_{2}\overline{u_{1}}^{6}\overline{u_{2}}^{5}), \\
\pi_{\tau_{6}}^{*}f(u_{1},u_{2}) &=& u_{1}^{8}u_{2}^{4}\overline{u_{1}}^{4}\overline{u_{2}}^{2}(1+\underline{u_{1}^{2a+3b-8}u_{2}^{a+2b-4}\overline{u_{1}}^{13-2a-3b}\overline{u_{2}}^{8-a-2b}}+u_{1}u_{2}^{2}+u_{1}u_{2}^{2}\overline{u_{1}}^{5}\overline{u_{2}}^{4}).
\end{eqnarray*}

For the exponents of the underlined terms, we calculate 
$$
\begin{array}{ll}
(2a+b-3)+(8-2a-b) = 5, \        &(3a+2b-6)+(12-3a-2b) = 6, \\
(3a+2b-6)+(12-3a-2b) = 6, \   &(a+b-3)+(5-a-b) = 2, \\
(a+b-3)+(5-a-b) = 2, \            &(3a+4b-12)+(18-3a-4b) = 6, \\
(3a+2b-12)+(18-3a-4b) = 6, \  &(2a+3b-8)+(13-2a-3b) = 5, \\
(2a+3b-8)+(13-2a-3b) = 5, \    &(a+2b-4)+(8-a-2b) = 4.
\end{array}
$$

Hence, we define the value of each underlined term as follows: 
If the exponent of $u_1$ or $\overline{u_1}$ is negative, %%
then we define the value of the term at $(0, u_2)$ to be $0$. 
If the exponent of $u_2$ or $\overline{u_2}$ is negative, %%
then we define the value of the term at $(u_1, 0)$ to be $0$. 

Then we see that all underlined terms are of class $C^1$. 
Hence, $\tilde{R}$ and $\tilde{f}$ are also of class $C^{1}$. 

However, by the above calculations, 
we have $\tilde{f}(0,0)\neq 0$, and we conclude that 
$\tilde V(\tau)^* = \emptyset$ in some neighborhood of $(u_1, u_2) = (0,0)$. %%%%

By the above results, %%%%
we see that $f$ and $\Sigma^*$ satisfy {\sc Assumption}{\rm (*)} in the statement of Theorem \ref{theorem11-improved}. 
Hence, by Theorem \ref{theorem11-improved}, 
The strict transform $\tilde V$ is 
a $C^1$-manifold in some neighborhood of $\hat{\pi}^{-1}(\zero) \cap \tilde V$ 
for every $(a,b)$ with $0 \leq a \leq 4$ and $0\leq b \leq 3$ 
and 
a real analytic manifold outside of $\hat{\pi}^{-1}(\zero) \cap \tilde V$. 

\bigskip

\bigskip

We now summarize the contents of this section as follows: 

\medskip

{\em 
Each of mixed polynomial germs 
$$
f_{a,b}(\bbz, \bar{\bbz}) := z_1^4 \bar{z}_1^2 + z_1^a \bar{z}_1^{4-a} z_2^b \bar{z}_2^{3-b} + \bar{z}_1^2 z_2^3 + z_2^3 \bar{z}_2^3
$$
at $\zero \in \bc^2$, where $0 \leq a \leq 4, \ 0\leq b \leq 3$ are integers, 
has five face functions. 
The four face functions 
$$
{(f_{a,b})}_U = {(f_{a,b})}_T = z_1^4\bar{z}_1^2,\ \ {(f_{a,b})}_S = \bar{z}_1^2z_2^3,\ \ {(f_{a,b})}_Q = z_1^4\bar{z}_1^2 + \bar{z}_1^2z_2^3,\ \ 
{(f_{a,b})}_P = \bar{z}_1^2z_2^3 + z_2^3\bar{z}_2^3
$$
are strongly polar non-negative mixed weighted homogeneous polynomials. 
Especially, the two face functions 
${(f_{a,b})}_P = \bar{z}_1^2z_2^3 + z_2^3\bar{z}_2^3$ and ${(f_{a,b})}_R = z_2^3\bar{z}_2^3$ 
are strongly mixed weighted homogeneous polynomials of polar degree $0$. 
Moreover, 
$(f_{a,b},\zero)$ is strongly Newton non-degenerate over $\De(U)=\De(T),\ \De(S),\ \De(Q)$ and $\De(P)$, 
however, it is not strongly Newton non-degenerate but Newton non-degenerate over $\De(R)$. 

Let $V := f_{a,b}^{-1}(0)$ be the germ of the mixed hypersurface at $\zero \in \bc^2$ and 
$\tilde V$ be the strict transform of $V$ to $X$ via the canonical toric modification \eqref{modification_exm4_3}. 
Then, $(f_{a,b},\zero)$, $\Sigma^*$ and $\tilde V$ satisfy the assumptions of Theorem \ref{theorem11-improved}. 
By Theorem \ref{theorem11-improved}, 
the strict transform $\tilde V$ is 
a \underline{$C^1$-manifold} in some neighborhood of $\hat{\pi}^{-1}(\zero) \cap \tilde V$ 
and a real analytic manifold outside of $\hat{\pi}^{-1}(\zero) \cap \tilde V$ 
for every $(a,b)$ with $0 \leq a \leq 4$ and $0\leq b \leq 3$. 
}

\vspace{1.2cm}

\noindent
{\em Acknowledgments.} \hfill

The authors would like to thank 
the referee for careful readings and many good advices, 
Professor Mutsuo Oka (Tokyo University of Science) for his kind and good advices, 
and 
Professor Toru Ohmoto (Hokkaido University) for introducing them studies of mixed functions. 

\vspace{1.2cm}

%%%%%%%%%%%%%%%%%%%%%%%%%%%%%%%%%%%%%%%

\end{document}